\documentclass{amsart}
\usepackage[dvips]{epsfig}
\newtheorem{thm}{Theorem}[subsection]
 \newtheorem{cor}[thm]{Corollary}
 \newtheorem{lem}[thm]{Lemma}
 
 \theoremstyle{definition}
\newtheorem{fact}{Fact}
 \newtheorem{formula}{Formula}
\newtheorem{defn}[thm]{Definition}
 \newtheorem{rem}{Remark}
 \numberwithin{equation}{subsection}

\begin{document}
\title[Positive Dehn Twist Expressions for some New Involutions]
 {Positive Dehn Twist Expressions for some New Involutions in Mapping Class Group}

\author{ Yusuf Z. Gurtas }

\address{Department of Mathematics, Suffolk CCC, Selden, NY, USA}

\email{gurtasy@sunysuffolk.edu}

\thanks{}




\dedicatory{}

\commby{}

\begin{abstract}
The well-known fact that any genus $g$ symplectic Lefschetz fibration $%
X^{4}\rightarrow S^{2}$ is given by a word that is equal to the
identity element in the mapping class group and each of whose
elements is given by a positive Dehn twist, provides an intimate
relationship between words in the mapping class group and
4-manifolds that are realized as symplectic Lefschetz fibrations.
In this article we provide new words in the mapping class group,
hence new symplectic Lefschetz fibrations. We also compute the
signatures of those symplectic Lefschetz fibrations.

\end{abstract}

\maketitle

\section*{Introduction}

The last decade has experienced a resurgence of interest in the
mapping class groups of 2-dimensional closed oriented surfaces.
This is primarily
due to Donaldson's theorem \cite{Do} that any closed oriented symplectic $4-$%
manifold has the structure of a Lefschetz pencil and Gompf's
theorem \cite{GS} that any Lefschetz pencil supports a symplectic
structure. These theorems, together with the fact that any genus
$g$ symplectic Lefschetz fibration is given by a word that is
equal to the identity element in the mapping class group and each
of whose elements is given by a positive Dehn twist, provides an
intimate relationship between words in the mapping class group and $4-$%
dimensional symplectic topology. Therefore the elements of finite
order in the mapping class group are of special importance. There
are very few examples of elements of finite order for which
explicit positive Dehn twist products are known and a great deal
is known about the structure of the $4-$manifolds that they
describe. In this article we give an algorithm for positive Dehn
twist products for a set of new involutions in the mapping class
group that are non-hyperelliptic. These involutions are obtained
by combining the positive Dehn twist expressions for two
well-known involutions of the mapping class group.

An application of the words in the mapping class group that we
produced is to determine the homeomorphism types of the
$4-$manifolds that they describe. This calculation involves
computing two invariants of those manifolds. The first one, which
is very easy to compute, is the Euler characteristic, and the
second one is the signature. Using the algorithm described in
\cite{Oz} we wrote a Matlab program that computes the signatures
of some of these manifolds. The computations that we have done
using this program point to a closed formula for the signature,
which is mentioned in the last part of this article.

This article will be followed by two other articles, the first one
of which will contain the same computations for the multiple case.
Namely, the explicit positive Dehn twist products for the
involutions obtained by combining several involutions of the type
described in this article and the signature computations for the
4-manifolds they describe as symplectic Lefschetz fibrations will
be presented. The third article will contain the same work for
some finite order elements in the mapping class group that are
rotations through $2\pi/p, p-$odd.

\section{Mapping Class Groups}

\subsection{Definition of $M_g$}

Let $\Sigma_{g}$ be a $2$-dimensional, closed, compact, oriented
surface of genus $g>0$. From now on, $g$ will always be positive
unless otherwise stated. Let $Diff^{+}\left( \Sigma_{g}\right) $
be the group of all orientation-preserving diffeomorphisms
$\Sigma_{g}\rightarrow \Sigma_{g},$ and $ Diff_{0}^{+}\left(
\Sigma_{g}\right) $ be the subgroup of $Diff^{+}\left(
\Sigma_{g}\right) $ consisting of all orientation-preserving diffeomorphisms $%
\Sigma_{g}\rightarrow \Sigma_{g}$ that are isotopic to the
identity.

The mapping class group $M_{g}$ of $\Sigma_{g}$ is defined to be
the group of isotopy classes of orientation-preserving
diffeomorphisms of $\Sigma_{g}$, i.e.,
\[
M_{g}=Diff^{+}\left( \Sigma_{g}\right) /Diff_{0}^{+}\left(
\Sigma_{g}\right) .
\]

\subsection{Generators}

\begin{defn}
Let $\alpha $ be a simple closed curve on $\Sigma_{g}$. Cut the
surface $\Sigma_{g}$ open along $\alpha $ and glue the ends back
after
rotating one of the ends $360^{\circ }$ to the right, which makes sense if $%
\Sigma_{g}$ is oriented. This operation defines a diffeomorphism $%
\Sigma_{g}\rightarrow \Sigma_{g}$ supported in a small
neighborhood of $\alpha $ and \emph{the positive Dehn twist about
}$\alpha $ is defined to be the isotopy class of this
diffeomorphism. It is denoted by $t_{\alpha }$ or $ D\left( \alpha
\right) .$
\end{defn}

\begin{fact}
$t_{\alpha }$ is independent of how $\alpha$ is oriented.
\end{fact}

A simple closed curve $\alpha $ on $\Sigma_{g}$ is called
\emph{nonseparating} if $\Sigma_{g}\backslash \alpha $ has one
connected component and it is called \emph{separating} otherwise.

\begin{fact}
 $M_{g}$ is generated by Dehn twists about nonseparating simple
closed curves.
\end{fact}

This fact was first proven by Dehn in 1938 \cite{De}. His set of
generators consisted of $2g\left( g-1\right) $ generators. Much
later (1964) Lickorish independently proved that $3g-1$ Dehn
twists are enough to generate $M_{g}$ \cite{Li}. Humphries (1979)
proved that $2n+1$ Dehn twists are enough and this is the minimal
number of twists to generate $M_{g}$ \cite{Hu}. If we do not
require the generators to be Dehn twists then the mapping class
groups can be generated by $2$ elements \cite{W2}.

This is clearly the lowest number of generators as we know that
$M_{g}$ is not an abelian group.

\subsection{Presentation} First, we need a few technical lemmas
that will be useful in proving some of the defining relations in
the presentation for $M_g$.

\begin{lem} \label{f(a)=b}
Let $\alpha $ and $\beta $ be two non-separating simple
closed curves on $\Sigma_{g}.$ Then there is a diffeomorphism $%
f:\Sigma_{g}\rightarrow \Sigma_{g}$ such that $f\left( \alpha
\right) =\beta .$
\end{lem}

For a proof see \cite{I}.

Let $ \mathcal{S} \left( \Sigma_{g}\right) $ be the set of all
isotopy classes of simple closed curves in $\Sigma_{g}.
$ For $[\alpha ],$ $%
[\beta ]$ $\in \mathcal{S}\left( \Sigma_{g}\right) $ define
$I\left( [\alpha ],[\beta ]\right) =\min \{|a\cap b|$ $|a\in
[\alpha ],b\in [\beta ]\}.$

\begin{defn}
If two simple closed curves $\alpha $ and $\beta $ intersect each
other transversely at one point, denote it by $\alpha \bot $
$\beta ,$ then we define the product $\alpha \beta $ as the simple
closed curve obtained from $\alpha \cup \beta $ by resolving the
intersection point $\alpha \cap \beta $ according to Figure
\ref{resolve}.
 \begin{figure}[htbp]
     \centering  \leavevmode
     \psfig{file=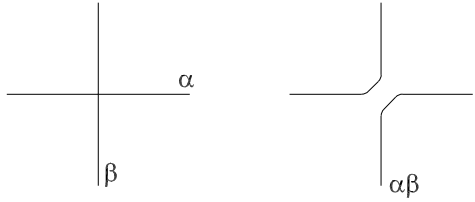,width=2.50in,clip=}
     \caption{Resolving an Intersection}
     \label{resolve}
 \end{figure}

\end{defn}

If $|\alpha \cap \beta |=I\left( [\alpha ],[\beta ]\right) =k$
then $\alpha \beta $ is defined as the simple closed curve
obtained from $\alpha \cup \beta $ by resolving each of the $k$
intersection points according to Figure \ref{resolve}.

\noindent {\bf Notation:} \thinspace   $\alpha (\alpha \beta
)=\alpha ^{2}\beta ,$ $\alpha (\alpha (\alpha \beta ))=\alpha
^{3}(\beta ),\ldots ,\alpha (\alpha \cdots (\alpha \beta ))=\alpha
^{k}\beta .$

\begin{rem}
\begin{enumerate}
\item[$\left( a\right) $]  If $I\left( [\alpha ],[\beta ]\right) =1,$ then $%
t_{\alpha }\left( \beta \right) =\alpha \beta $ and $t_{\beta
}\left( \alpha \right) =\beta \alpha .$ \label{a(b)=ab}

\item[$\left( b\right) $]  If $I\left( [\alpha ],[\beta ]\right) =0,$ then $%
t_{\alpha }\left( \beta \right) =\beta $ and $t_{\beta }\left(
\alpha \right) =\alpha .$
\end{enumerate}
\end{rem}

\begin{lem} \label{intnum.lem}
\begin{enumerate}
\item[$\left( a\right) $]  If $\alpha $ and $\beta $ are two
simple closed curves, then $t_{\alpha }\left( \beta \right)
=\alpha ^{k}\beta ,$ where $ I\left( [\alpha ],[\beta ]\right)
=k,$ (\cite{Lu2}).

\item[$\left( b\right) $]  $I\left( [\alpha ],[\beta ]\right)
=I\left( [\alpha \beta ],[\beta ]\right) =I\left( [\alpha \beta
],[\alpha ]\right) .$
\end{enumerate}
\end{lem}

\begin{lem} \label{a(ba)=b.lem}
If $I\left( [\alpha ],[\beta ]\right) =1,$ then $\alpha \left(
\beta \alpha \right) =\beta $ and $\beta \left( \alpha \beta
\right)
=\alpha ,$ namely, $t_{\alpha }t_{\beta }\left( \alpha \right) =\beta $ and $%
t_{\beta }t_{\alpha }\left( \beta \right) =\alpha .$
\end{lem}

\noindent {\bf Proof:} The following figure is a sketch of the
proof.

 \begin{figure}[htbp]
     \centering  \leavevmode
     \psfig{file=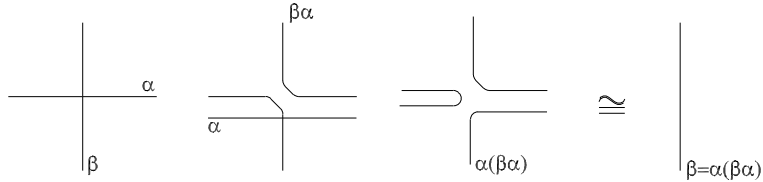,width=4.0in,clip=}
     \caption{$t_{\alpha }t_{\beta }\left( \alpha \right) =\beta $ and $%
t_{\beta }t_{\alpha }\left( \beta \right) =\alpha $, for $I\left(
[\alpha ],[\beta ]\right) =1$}
     \label{braided.fig}
 \end{figure}

\begin{lem} \label{conj.lem}
For any simple closed curve $\alpha $ in $%
\Sigma_{g}$ and any diffeomorphism $f:\Sigma_{g}\rightarrow
\Sigma_{g},$ we have $t_{f\left( \alpha \right) }=ft_{\alpha
}f^{-1}.$
\end{lem}

For a proof see \cite{I}.

\begin{cor} \label{conj.cor}
The Dehn twists about any two nonseparating simple closed curves
$\alpha $ and $\beta $ are conjugate, i.e., $t_{\beta }=ft_{\alpha
}f^{-1}$, for some $f:\Sigma_{g}\rightarrow \Sigma_{g}$.
\end{cor}

This easily follows from the previous lemma and Lemma
\ref{f(a)=b}.

\begin{rem} \label{conj.rem}
\begin{enumerate}
\item[$\left( a\right) $]  In particular if $f=t_{\beta }$ in
Lemma \ref{conj.lem}, then $t_{t_{\beta }\left( \alpha \right)
}=t_{\beta }t_{\alpha }t_{\beta }^{-1}.$

\item[$\left( b\right) $]  Moreover, if $I\left( [\alpha ],[\beta
]\right)
=1,$ then $t_{\beta \alpha }=t_{\beta }t_{\alpha }t_{\beta }^{-1},$ because $%
t_{\beta }\left( \alpha \right) =\beta \alpha $ by Remark
\ref{a(b)=ab} (a) .
\end{enumerate}
\end{rem}

\noindent {\bf Commutativity and Braid Relations}

\begin{lem} \label{com&braid.lem}
\begin{enumerate}
\item[$\left( a\right) $]  If $I\left( \left[ \alpha \right]
,\left[ \beta \right] \right) =0,$ then $t_{\alpha }t_{\beta
}=t_{\beta }t_{\alpha }.$

\item[$\left( b\right) $]  If $I\left( [\alpha ],[\beta ]\right) =1,$ then $%
t_{\alpha }t_{\beta }t_{\alpha }=t_{\beta }t_{\alpha }t_{\beta }.$
\end{enumerate}
\end{lem}

\noindent {\bf Proof:} \thinspace
\begin{enumerate}
\item [$\left( a\right) $] We have $t_{\beta }\left( \alpha
\right) =\alpha $ by Remark \ref{a(b)=ab} (b) . Therefore,
$t_{t_{\beta }\left( \alpha \right) }=t_{\beta }t_{\alpha
}t_{\beta }^{-1}$ in Remark \ref{conj.rem} (a) becomes $t_{\alpha
}=t_{\beta }t_{\alpha }t_{\beta }^{-1}.$

\item[$\left( b\right) $]  Since $I\left( [\alpha ],[\beta
]\right) =1$ we have $I\left( [\alpha ],[\beta \alpha ]\right) =1$
by Lemma \ref{intnum.lem} (b). Therefore, we also have $t_{\alpha
}\left( \beta \alpha \right) =\alpha \left( \beta \alpha \right) $
using Remark \ref{a(b)=ab} (a), which is equal to $\beta $, using
Lemma \ref{a(ba)=b.lem}. Now, the twist along $\beta $ can be
expressed as
\[
t_{\beta }=t_{t_{\alpha }\left( \beta \alpha \right) }=t_{\alpha
}t_{\beta \alpha }t_{\alpha }^{-1},
\]
thanks to Lemma \ref{conj.lem}. Finally substituting $t_{\beta
}t_{\alpha }t_{\beta }^{-1}$ for $t_{\beta \alpha }$ in the last
equation, using Remark \ref {conj.rem} (b), we get
\[
t_{\beta }=t_{\alpha }t_{\beta }t_{\alpha }t_{\beta
}^{-1}t_{\alpha }^{-1},
\]
which finishes the proof.
\end{enumerate}

The relation $t_{\alpha }t_{\beta }t_{\alpha }=t_{\beta }t_{\alpha
}t_{\beta}$, for $\left[ \alpha \right] \neq \left[ \beta \right]
$, is called the \emph{braid relation} and we say $t_{\alpha }$ and $%
t_{\beta }$ are \emph{braided} if they satisfy the braid relation.

The following formula plays the key role in proving the converse
of Lemma \ref{com&braid.lem}, (\cite{FLP}).

\begin{formula} \label{formula1}
For two simple closed curves $\alpha $ and $\beta $ we have
\[
I\left( [t_{\alpha }^{n}\left( \beta \right) ],[\beta ]\right)
=\left| n\right| I\left( \left[ \alpha \right] ,\left[ \beta
\right] \right) ^{2}.
\]
\end{formula}

In addition to Formula \ref{formula1}, the following two simple
facts will also be used in proving the converse of Lemma
\ref{com&braid.lem}.

\begin{fact} \label{isotopcurves}
For two simple closed curves $\alpha $ and $\beta,$ if $\left[
\alpha \right] \neq \left[ \beta \right]$, then $ t_{\alpha }\neq
t_{\beta },$ (\cite{PR}).
\end{fact}

\begin{fact} \label{intnonzero}
For two simple closed curves $\alpha $ and $\beta $ if $%
I\left( \left[ \alpha \right] ,\left[ \beta \right] \right) \neq 0,$ then $%
\left[ t_{\alpha }\left( \beta \right) \right] \neq \left[ \beta
\right] $.
\end{fact}

\noindent {\bf Proof:} \thinspace Formula \ref{formula1} for $n=1$
gives $I\left( [t_{\alpha }\left( \beta \right) ],[\beta ]\right)
=I\left( \left[ \alpha \right] ,\left[ \beta \right] \right)
^{2}\neq 0.$ If $ \left[ t_{\alpha }\left( \beta \right) \right] =
\left[ \beta \right] $ then $ I\left( [t_{\alpha }\left( \beta
\right) ],[\beta ]\right) =I\left( [\beta ],[\beta ]\right) $.
Since $I\left( [\beta ],[\beta ]\right) =0$ we can't have $\left[
t_{\alpha }\left( \beta \right) \right] =\left[ \beta \right].$

\begin{lem}
\begin{enumerate}
\item[$\left( a\right) $]  For two simple closed curves $\alpha $
and $\beta $ if  $t_{\alpha }t_{\beta }=t_{\beta }t_{\alpha },$
then $I\left( \left[ \alpha \right] ,\left[ \beta \right] \right)
=0.$

\item[$\left( b\right) $]  For two simple closed curves $\alpha $
and $\beta, $ with $\left[ \alpha \right] \neq \left[ \beta
\right], $ if $t_{\alpha }t_{\beta }t_{\alpha }=t_{\beta
}t_{\alpha }t_{\beta }$ then $I\left( \left[ \alpha \right]
,\left[ \beta \right] \right) =1.$
\end{enumerate}
\end{lem}

\noindent {\bf Proof:} \thinspace
\begin{enumerate}
\item[$\left( a\right) $]  From $t_{\alpha }t_{\beta }=t_{\beta
}t_{\alpha }$ we get $t_{\alpha }t_{\beta }t_{\alpha
}^{-1}=t_{\beta }.$ Using Lemma \ref{conj.lem}
we conclude that $t_{t_{\alpha }\left( \beta \right) }=t_{\beta }.$ So, $%
\left[ t_{\alpha }\left( \beta \right) \right] =\left[ \beta
\right] $ by Fact \ref{isotopcurves} and therefore $I\left( \left[
\alpha \right] ,\left[ \beta \right] \right) =0$ using Fact
\ref{intnonzero}.

\item[$\left( b\right) $]  From $t_{\alpha }t_{\beta }t_{\alpha
}=t_{\beta }t_{\alpha }t_{\beta }$ we get $t_{\alpha }t_{\beta
}t_{\alpha }\left(
t_{\alpha }t_{\beta }\right) ^{-1}=t_{\beta }.$ Using Lemma \ref{conj.lem} we get $%
t_{t_{\alpha }t_{\beta }\left( \alpha \right) }=t_{\beta }$, and therefore $%
\left[ t_{\alpha }t_{\beta }\left( \alpha \right) \right] =\left[
\beta \right] $ by Fact \ref{isotopcurves}. Applying Formula
\ref{formula1} with $n=1$ we get
\[
I\left( \left[ \alpha \right] ,\left[ \beta \right] \right)
^{2}=I\left( [t_{\beta }\left( \alpha \right) ],[\alpha ]\right)
=I\left( [t_{\alpha }t_{\beta }\left( \alpha \right) ],[\alpha
]\right) =I\left( \left[ \beta \right] ,\left[ \alpha \right]
\right) =I\left( \left[ \alpha \right] ,\left[ \beta \right]
\right) .
\]
 Therefore $I\left( \left[
\alpha \right] ,\left[ \beta \right] \right) $ is $0$ or $1.$ If
it is $0,$ then $t_{\alpha }$ and $t_{\beta }$ commute by Lemma
\ref {com&braid.lem} and $t_{\alpha }t_{\beta }t_{\alpha
}=t_{\beta }t_{\alpha }t_{\beta }$ becomes $t_{\beta }t_{\alpha
}^{2}=t_{\beta }^{2}t_{\alpha }$, and therefore $t_{\alpha
}=t_{\beta },$ i.e., $\left[ \alpha \right] =\left[ \beta \right]
,$ which is a contradiction.
\end{enumerate}

For more results on the previous two Lemmas see \cite{Mar}.

\noindent {\bf Lantern and Chain Relations}

The following two relations will be used for the presentation of
$M_{g},$ along with the commutativity and braid relations.

\begin{figure}[htbp]
     \centering  \leavevmode
     \psfig{file=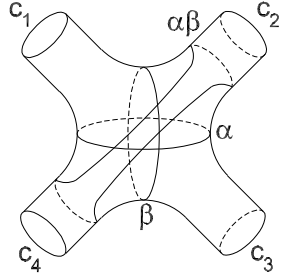,width=4.0in,clip=}
     \caption{The Lantern Relation}
     \label{lantern.fig}
 \end{figure}

Let $\Sigma_{0,4}$ be a sphere with four holes. If $c_{1},c_{2},c_{3},$ and $%
c_{4}$ are the boundary curves of $\Sigma_{0,4}$ and the simple closed curves $%
\alpha $ and $\beta $ are as shown in Figure \ref{lantern.fig} ,
then we have
\[
t_{\alpha }t_{\beta }t_{\alpha \beta
}=t_{c_{1}}t_{c_{2}}t_{c_{3}}t_{c_{4}},
\]
where $t_{c_{i}},$ $1\leq i\leq 4,$ denote the Dehn twists about
$c_{i}.$

This relation was known to Dehn and later on was rediscovered by
D.Johnson and named as \emph{lantern relation} by him
\cite{I},\cite{Lu1},\cite{J}. For more results on lantern relation
see \cite{Mar}.

The following relation was also known to Dehn and it is called the
\emph{\ chain relation. }

\begin{figure}[htbp]
     \centering  \leavevmode
     \psfig{file=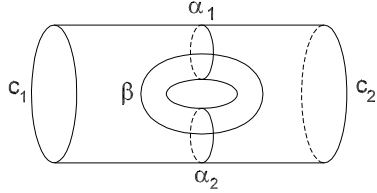,width=3.50in,clip=}
     \caption{The Chain Relation}
     \label{chain.fig}
 \end{figure}

Let $\Sigma_{1,2}$ be a torus with two boundary components. If $c_{1}$ and $%
c_{2}$ are the boundary curves of $\Sigma_{1,2}$ and $\alpha
_{1},$ $\alpha _{2},$ and $\beta $ are the simple closed curves as
shown in Figure \ref{chain.fig}, then we have

\[
\left( t_{\alpha _{1}}t_{\beta }t_{\alpha _{2}}\right)
^{4}=t_{c_{1}}t_{c_{2}}.
\]

\begin{figure}[htbp]
     \centering  \leavevmode
     \psfig{file=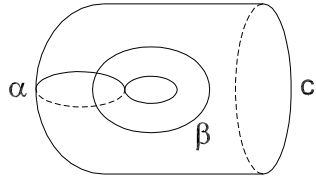,width=3.50in,clip=}
     \caption{Special Case of the Chain Relation}
     \label{specialchain.fig}
 \end{figure}

If $c _{1}$ bounds a disk then the chain relation becomes
\[
\left( t_{\alpha }t_{\beta }t_{\alpha }\right) ^{4}=t_{c},
\]
where $\alpha =\alpha _{1}=\alpha _{2}$ and $c=c_{2}$ is the only
boundary curve. This is a special case of the chain relation and
can be rewritten as
\[
\left( t_{\alpha }t_{\beta }\right) ^{6}=t_{c}
\]
using the braid relation $t_{\alpha }t_{\beta }t_{\alpha
}=t_{\beta }t_{\alpha }t_{\beta }$ twice, Figure
\ref{specialchain.fig}.

One can also show that $\left( t_{\beta }t_{\alpha }\right)
^{6}=t_{c},$ using the braid relation.

If we call an ordered set of $n$ simple closed curves $c_n$, where
$c_i\bot c_{i+1}$ and $c_i \cap c_j=\emptyset $ for $|i-j|>1$,
\emph{a chain of length n} \cite{W3}, then the chain relation
defined above is realized by a chain of length 3 and the special
case of that relation shown in Figure \ref{specialchain.fig} is
realized by a chain of length 2.

 A simple and explicit presentation for $M_{g}$, $g\geq
3,$ was not given until 1983 \cite{W1}. Wajnryb gave the following
presentation following the ideas of Hatcher-Thurston and Harer.

\begin{figure}[htbp]
     \centering  \leavevmode
     \psfig{file=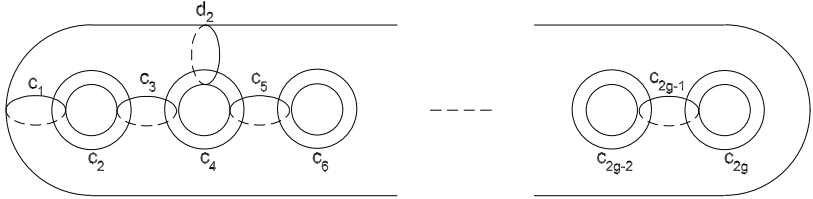,width=4.50in,clip=}
     \caption{The Generators}
     \label{gener.fig}
 \end{figure}

\begin{figure}[htbp]
     \centering  \leavevmode
     \psfig{file=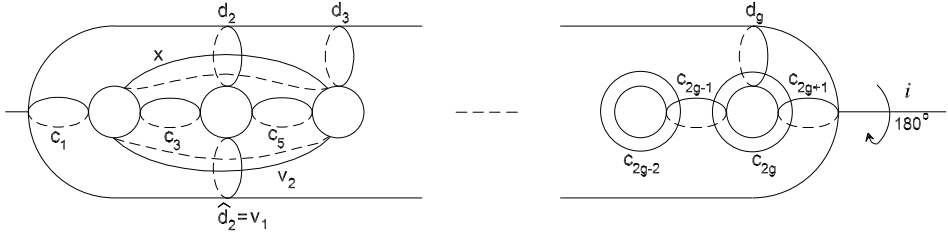,width=5.0in,clip=}
     \caption{The Relations and the Hyperelliptic Involution}
     \label{relat.fig}
 \end{figure}

\begin{thm} \label{present.thm}
The mapping class group $M_{g}$ of a $2$-dimensional, closed,
compact, oriented surface $\Sigma_{g}$ of genus $g\geq 3$ admits a
presentation with generators
\[
t_{d_{2}},t_{c_{1}},t_{c_{2}},\ldots ,t_{c_{2g}},
\]
where the cycles $d_{2},c_{1},c_{2},\ldots ,c_{2g}$ are as shown
in Figure \ref{gener.fig}, and with the following defining
relations:

\begin{enumerate}
\item[$\left( A\right) $]  $t_{d_{2}}$ and $t_{c_{4}}$ are braided and $%
t_{d_{2}}$ commutes with $t_{c_{i}}$ for $i\neq 4.$ $t_{c_{i}}$
commutes with $t_{c_{j}}$ if $\left| i-j\right| >1,$ and $t_{c_{i}}$ and $%
t_{c_{j}}$ are braided if $\left| i-j\right| =1,$ $1\leq i,j\leq
2g,$ Figure \ref{gener.fig}.

\item[$\left( B\right) $]  $t_{c_{1}},t_{c_{2}},t_{c_{3}},t_{d_{2}}$ and $t_{%
\widehat{d}_{2}}$ satisfy the chain relation, i.e., $\left(
t_{c_{1}}t_{c_{2}}t_{c_{3}}\right)
^{4}=t_{d_{2}}t_{\widehat{d}_{2}},$ where $\widehat{d}
_{2}=t_{c_{4}}t_{c_{3}}t_{c_{2}}t_{c_{1}}^{2}t_{c_{2}}t_{c_{3}}t_{c_{4}}%
\left( d_{2}\right) ,$ Figure \ref{gener.fig},\ref{relat.fig}.

\item[$\left( C\right) $]  $%
t_{c_{1}},t_{c_{3}},t_{c_{5}},t_{d_{3}},t_{d_{2}},t_{x},$ and
$t_{d_{2}x}$ satisfy the lantern relation, i.e.,
\[
t_{c_{1}}t_{c_{3}}t_{c_{5}}t_{d_{3}}=t_{d_{2}}t_{x}t_{d_{2}x},
\]
where $x=t_{r_{1}}t_{r_{2}}\left( d_{2}\right) ,$ $
t_{r_{1}}=t_{c_{2}}t_{c_{1}}t_{c_{3}}t_{c_{2}},$ $%
t_{r_{2}}=t_{c_{4}}t_{c_{3}}t_{c_{5}}t_{c_{4}},$ $d_{3}=$ $%
t_{c_{6}}t_{c_{5}}t_{c_{4}}t_{c_{3}}t_{c_{2}}t_{a}\left( b\right) ,$ $%
a=(t_{c_{4}}t_{c_{3}}t_{c_{5}}t_{c_{4}}t_{c_{6}}t_{c_{5}})^{-1}\left(
d_{2}\right) ,$ and $b=\left(
t_{c_{4}}t_{c_{3}}t_{c_{2}}t_{c_{1}}\right) ^{-1}\left(
d_{2}\right) $, Figure \ref{relat.fig}.

\item[$\left( D\right) $] $t_{c_{1}},t_{c_{2}},\ldots,t_{c_{2g}}$
and $t_{d_{g}}$ satisfy the relation
\[
[t_{c_{2g}}t_{c_{2g-1}}\cdots
t_{c_{2}}t_{c_{1}}^{2}t_{c_{2}}\cdots
t_{c_{2g-1}}t_{c_{2g}},t_{d_{g}}]=1,
\]
where $d_{g}=t_{u_{g-1}}t_{u_{g-2}}\cdots t_{u_{1}}\left( c_{1}\right) ,$ $%
t_{u_{i}}=\left(
t_{c_{2i-1}}t_{c_{2i}}t_{c_{2i+1}}t_{c_{2i+2}}\right)
^{-1}t_{v_{i}}t_{c_{2i+2}}t_{c_{_{2i+1}}}t_{c_{_{2i}}}$ for $1\leq
i\leq g-1, $ $v_{1}=\widehat{d}_{2},$
$v_{i}=t_{r_{i-1}}^{-1}t_{r_{i}}^{-1}\left(
v_{i-1}\right) $ for $2\leq i\leq g-1,$ $%
t_{r_{i}}=t_{c_{_{2i}}}t_{c_{_{2i+1}}}t_{c_{_{2i-1}}}t_{c_{_{2i}}}$ for $1%
\leq i\leq g-1,$ Figure \ref{relat.fig}.
\end{enumerate}
\end{thm}

Relation $\left( D\right) $ is related to the so called
\emph{hyperelliptic involution} $i:\Sigma_{g}\rightarrow
\Sigma_{g}$ of $\Sigma_{g}.$ It is the element of order $2,$
geometrically represented as the $180^{\circ }$ rotation about the
horizontal axis as shown in Figure \ref{relat.fig}. It fixes the
unoriented cycles $c_{1},c_{2},\ldots ,c_{2g},d_{g}=$ $c_{2g+1}$
and acts as $-Id$ on the homology. It can be shown that $i$ can be
expressed as $i=t_{c_{2g+1}}t_{c_{2g}}t_{c_{2g-1}}\cdots
t_{c_{2}}t_{c_{1}}^{2}t_{c_{2}}\cdots
t_{c_{2g-1}}t_{c_{2g}}t_{c_{2g+1}}.$


From the expression for $i$ it is clear that $i$ commutes with
$t_{d_{g}}$ if and only if $\left( D\right) $ holds.

\subsection{Involutions} \label{inv.sec}

The hyperelliptic involution $i:\Sigma_{g}\rightarrow \Sigma_{g}$
that is depicted in Figure \ref{relat.fig} is probably the most
studied and the best understood order $2$ element in the
mapping class group. The so-called \emph{hyperelliptic mapping class group} $%
H_{g}$, is the infinite subgroup of $M_{g}$ consisting of the
elements that commute with $i$. All of the generators of
$M_{g}$ commute with $i$ for $%
g=1,2,$ therefore $M_{g}=H_{g},$ for $g=1,2$.

For higher genus, $H_{g}$ is generated by $t_{c_{i}},$ $i=1,\ldots
,2g+1,$ where $c_{i}$ are the cycles that are shown in Figure
\ref{relat.fig} . The relations are

\begin{align}
\nonumber t_{c_{i}}t_{c_{j}}=t_{c_{j}}t_{c_{i}}\, \mbox{if} \,
\left| i-j\right|
>1,\\
\nonumber t_{c_{i}}t_{c_{i+1}}t_{c_{i}} =t_{c_{i+1}}t_{c_{i}}t_{c_{i+1}}, \\
\nonumber \left( t_{c_{2g+1}}t_{c_{2g}}\cdots
t_{c_{2}}t_{c_{1}}^{2}t_{c_{2}}\cdots
t_{c_{2g}}t_{c_{2g+1}}\right) ^{2}=1, \\
\nonumber \mbox{and}
 \left(t_{c_{1}}t_{c_{2}}\cdots t_{c_{2g}}t_{c_{2g+1}}\right) ^{2g+2} =1.
\end{align}

In particular this gives a presentation for $M_{2}$ because
$M_{2}=H_{2}.$

\begin{figure}[htbp]
     \centering  \leavevmode
     \psfig{file=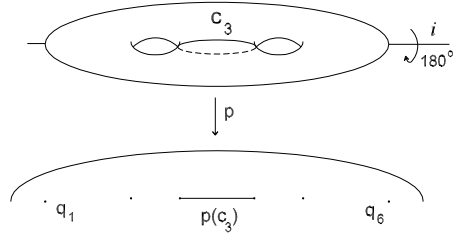,width=3.0in,clip=}
     \caption{Two Fold Branched Cover Defined by $i$}
     \label{proj.fig}
 \end{figure}

The orbit space of the involution $i$ is the sphere with $6$
marked points, denote it by $\Sigma_{0,6}.$ It defines a $2-$ fold
branched covering $p:\Sigma_{2}\rightarrow \Sigma_{0,6},$ branched
at $6$ points as shown in Figure \ref{proj.fig}. The cycles
$c_{i},i=1,\ldots ,5$ in $\Sigma_{2}$ project to the segments
$p\left( c_{i}\right) $ in $\Sigma_{0,6}$ that connect the marked
points $q_{i}$ and $q_{i+1}$. The isotopy classes of the half
twists
about $p\left( c_{i}\right) ,$ which we denote by $%
w_{i}$ for $i=1,\ldots ,5,$ generate the mapping class group
$M\left( S^{2},6\right) $ of $\Sigma_{0,6}$ (The marked points are
fixed setwise). There is a surjective homomorphism
\[
\psi :H_{2}\rightarrow M\left( S^{2},6\right)
\]
sending $t_{c_{i}}$ to $w_{i},i=1,\ldots ,5,$ and with $\ker \psi
=<i>.$ For a proof see \cite{B}.

\begin{figure}[htbp]
     \centering  \leavevmode
     \psfig{file=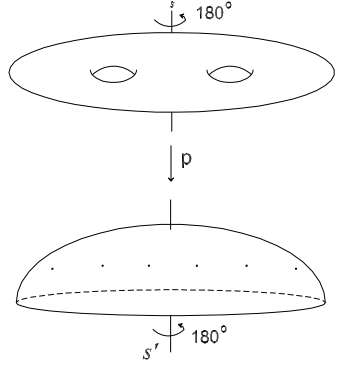,width=4.50in,clip=}
     \caption{The Involution $s$}
     \label{vertinv.fig}
 \end{figure}

There is another involution in the mapping class group which is
also geometric. We will denote it by $s.$ It is described as
$180{{}^{o}}$ rotation about the vertical axis as shown in Figure
\ref{vertinv.fig} . We will begin with finding an explicit
positive Dehn twist expression for $s\ $in $M_{2}.$ The following
paragraph explains the idea, which generalizes to $M_{g}$ with
some extra work. See \cite{McP}, \cite{Ko} for details.

Since $s$ commutes with $i$ in $M_{2},$ it descends to a map $\widetilde{s}%
:\Sigma_{0,6}\rightarrow \Sigma_{0,6}.$ We will obtain an
expression for $\widetilde{s}$ in terms of the generators $w_{i}$
and then lift that expression to $H_{2}=M_{2}$ via the surjection
$\psi .$

\begin{figure}[htbp]
     \centering  \leavevmode
     \psfig{file=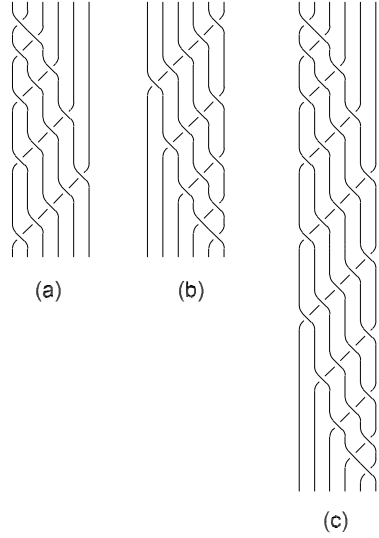, width=5.50in,clip=}
     \caption{Braids Representing $s'$ and $s'^{2}$ }
     \label{twobraids.fig}
 \end{figure}

$\widetilde{s}$ is counterclokwise rotation through $180{{}^{o}}$ on $%
\Sigma_{0,6},$ about the axis through the center and the south
pole of the sphere, fixing the marked points setwise. It fixes the
south pole,
therefore forgetting that point we can isotope $\widetilde{s}$ to a map $%
s^{\prime }$ of a disc including the $6$ marked points, Figure
\ref{vertinv.fig}. Being an element of $M\left( D,6\right)
=B_{6},$ $s^{\prime }$ can be realized as a braid. Figure
\ref{twobraids.fig} (a) and (b) show two different ways of
sketching that braid. The first one has the expression
\[
\sigma _{1}\sigma _{2}\sigma _{1}\sigma _{3}\sigma _{2}\sigma
_{1}\sigma _{4}\sigma _{3}\sigma _{2}\sigma _{1}\sigma _{5}\sigma
_{4}\sigma _{3}\sigma _{2}\sigma _{1}
\]
and the second one has the expression
\[
\sigma _{5}\sigma _{4}\sigma _{3}\sigma _{2}\sigma _{1}\sigma
_{5}\sigma _{4}\sigma _{3}\sigma _{2}\sigma _{5}\sigma _{4}\sigma
_{3}\sigma _{5}\sigma _{4}\sigma _{5}.
\]

Since Figure \ref{twobraids.fig} (a) and (b) are representing
isotopic braids, they correspond to the same element in $M\left(
D,6\right) .$ Therefore
\[
s^{\prime 2}=\sigma _{5}\sigma _{4}\sigma _{3}\sigma _{2}\sigma
_{1}\sigma _{5}\sigma _{4}\sigma _{3}\sigma _{2}\sigma _{5}\sigma
_{4}\sigma _{3}\sigma _{5}\sigma _{4}\sigma _{5}\underline{\sigma
_{1}}\underline{\sigma _{2}\sigma _{1}}\underline{\sigma
_{3}\sigma _{2}\sigma _{1}}\sigma _{4}\sigma _{3}\sigma _{2}\sigma
_{1}\sigma _{5}\sigma _{4}\sigma _{3}\sigma _{2}\sigma _{1}
\]
which is equal to
\[
\sigma _{5}\sigma _{4}\sigma _{3}\sigma _{2}\sigma _{1}\sigma
_{5}\sigma _{4}\sigma _{3}\sigma _{2}\underline{\sigma _{1}}\sigma
_{5}\sigma _{4}\sigma _{3}\underline{\sigma _{2}\sigma _{1}}\sigma
_{5}\sigma _{4}\underline{\sigma _{3}\sigma _{2}\sigma _{1}}\sigma
_{5}\sigma _{4}\sigma _{3}\sigma _{2}\sigma _{1}\sigma _{5}\sigma
_{4}\sigma _{3}\sigma _{2}\sigma _{1}.
\]
Thus
\[
s^{\prime 2}=\left( \sigma _{5}\sigma _{4}\sigma _{3}\sigma
_{2}\sigma _{1}\right) ^{6}
\]
using commutativity relations only. It is not difficult to see
that this is isotopic to a right handed twist about the boundary
of the disk $D$, Figure \ref{twobraids.fig} (c). Therefore we have

\begin{eqnarray*}
\widetilde{s}^{2} &=&\left( w_{5}w_{4}w_{3}w_{2}w_{1}\right) ^{6} \\
&=&1
\end{eqnarray*}
in $M\left( S^{2},6\right)$ and hence the lift $h$ of
$\widetilde{s}$ satisfies
\[
h^{2}=\left( t_{c_{5}}t_{c_{4}}t_{c_{3}}t_{c_{2}}t_{c_{1}}\right)
^{6},
\]
which is equal to $1$ in $M_{2}.$ One needs to check the action of
$s$ and
\[
h=t_{c_{1}}t_{c_{2}}t_{c_{1}}t_{c_{3}}t_{c_{2}}t_{c_{1}}t_{c_{4}}t_{c_{3}}t_{c_{2}}t_{c_{1}}t_{c_{5}}t_{c_{4}}t_{c_{3}}t_{c_{2}}t_{c_{1}}
\]
on the homology to see that the action of $s$ is the same as the
action of $h $, but not as that of $h\circ i$ \cite{McP}.
Therefore $s=h$ and
\[
s^{2}=\left( t_{c_{5}}t_{c_{4}}t_{c_{3}}t_{c_{2}}t_{c_{1}}\right)
^{6}.
\]

Using the braid and commutativity relations it is not difficult to
show that
\[
s^{2}=\left( t_{c_{1}}t_{c_{2}}t_{c_{3}}t_{c_{4}}t_{c_{5}}\right)
^{6}
\]
as well.

\begin{figure}[htbp]
     \centering  \leavevmode
     \psfig{file=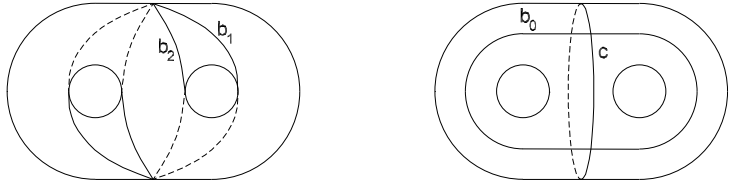, width=4.0in,clip=}
     \caption{Another Dehn Twist Expression for $s$ }
     \label{matsumotocurves}
 \end{figure}

 Another positive Dehn twist expression for $s,$ which is of
particular interest from topological point of view \cite{Ma}, is

\[
s=t_{b_{0}}t_{b_{1}}t_{b_{2}}t_{c},
\]
where $b_{0},b_{1},b_{2},$ and $c$ are the cycles shown in Figure
\ref{matsumotocurves}.

Distinguishing the cycle notation from the Dehn twist notation, it
is easy to see that
\begin{eqnarray*}
b_{0} &=&t_{c_{1}}t_{c_{2}}t_{c_{3}}t_{c_{4}}\left( c_{5}\right) , \\
b_{1} &=&t_{c_{1}}t_{c_{1}}t_{c_{2}}t_{c_{3}}\left( c_{4}\right) , \\
b_{2} &=&t_{c_{2}}t_{c_{1}}t_{c_{1}}t_{c_{2}}\left( c_{3}\right) .
\end{eqnarray*}
Applying Lemma \ref{conj.lem} we get
\begin{eqnarray*}
t_{b_{0}} &=&t_{c_{1}}t_{c_{2}}t_{c_{3}}t_{c_{4}}t_{c_{5}}\left(
t_{c_{1}}t_{c_{2}}t_{c_{3}}t_{c_{4}}\right) ^{-1}, \\
t_{b_{1}} &=&t_{c_{1}}t_{c_{1}}t_{c_{2}}t_{c_{3}}t_{c_{4}}\left(
t_{c_{1}}t_{c_{1}}t_{c_{2}}t_{c_{3}}\right) ^{-1}, \\
t_{b_{2}} &=&t_{c_{2}}t_{c_{1}}t_{c_{1}}t_{c_{2}}t_{c_{3}}\left(
t_{c_{2}}t_{c_{1}}t_{c_{1}}t_{c_{2}}\right) ^{-1}.
\end{eqnarray*}

 Therefore the product $t_{b_{0}}t_{b_{1}}t_{b_{2}}$
becomes
\[
t_{c_{1}}t_{c_{2}}t_{c_{3}}t_{c_{4}}t_{c_{5}}\left(
t_{c_{1}}t_{c_{2}}t_{c_{3}}t_{c_{4}}\right)
^{-1}t_{c_{1}}t_{c_{1}}t_{c_{2}}t_{c_{3}}t_{c_{4}}\left(
t_{c_{1}}t_{c_{1}}t_{c_{2}}t_{c_{3}}\right)
^{-1}t_{c_{2}}t_{c_{1}}t_{c_{1}}t_{c_{2}}t_{c_{3}}\left(
t_{c_{2}}t_{c_{1}}t_{c_{1}}t_{c_{2}}\right) ^{-1}.
\]

 We also have $t_{c}=\left( t_{c_{1}}t_{c_{2}}t_{c_{1}}\right)
^{4}=\left( t_{c_{1}}t_{c_{2}}\right) ^{6}=\left(
t_{c_{2}}t_{c_{1}}\right) ^{6}$ using the braid relation and the
special case of the chain relation on the subsurface $\Sigma_1^1$.

 Next, we will simplify the expression for the product
$t_{b_{0}}t_{b_{1}}t_{b_{2}}.$ We will use just the indices
representing the twists for brevity. For example $1$ will mean
$t_{c_{1}}$, $2$ will mean $t_{c_{2}}$, etc. The change in each
line occurs within the underlined portion of the entire expression
and the result of that change is not underlined in the next line.
We only use commutativity and braid relations.

\begin{eqnarray*}
t_{b_{0}}t_{b_{1}}t_{b_{2}} &=&12345\left( \underline{1}234\right) ^{-1}%
\underline{1}1234\left( 1123\right) ^{-1}21123\left( 2112\right) ^{-1} \\
&=&12345\underline{\left( 234\right) ^{-1}1234}\left( 1123\right)
^{-1}21123\left( 2112\right) ^{-1} \\
&=&123451234\left( 123\right) ^{-1}\left( \underline{11}23\right)
^{-1}21123\left( 2112\right) ^{-1} \\
&=&123451234\left( 123\right) ^{-1}\left( 23\right)
^{-1}\underline{\left(
11\right) ^{-1}}21123\left( 2112\right) ^{-1} \\
&=&123451234\left( 123\right) ^{-1}\left( \underline{2}3\right)
^{-1}21123\left( 11\right) ^{-1}\left( 2112\right) ^{-1} \\
&=&123451234\left( 123\right) ^{-1}\underline{\left( 3\right) ^{-1}}%
1123\left( 11\right) ^{-1}\left( 2112\right) ^{-1} \\
&=&123451234\left( \underline{1}23\right) ^{-1}1123\left( 2\right)
^{-1}\left( 11\right) ^{-1}\left( 2112\right) ^{-1} \\
&=&123451234\underline{\left( 23\right) ^{-1}}123\left( 2\right)
^{-1}\left(
11\right) ^{-1}\left( 2112\right) ^{-1} \\
&=&123451234123\left( 12\right) ^{-1}\left( 2\right) ^{-1}\left(
11\right)
^{-1}\left( 2112\right) ^{-1} \\
&=&123451234123\underline{1211^{-1}\left( 12\right) ^{-1}}\left(
12\right)
^{-1}\left( 2\right) ^{-1}\left( 11\right) ^{-1}\left( 2112\right) ^{-1} \\
&=&1234512341231211^{-1}2^{-1}1^{-1}2^{-1}1^{-1}2^{-1}1^{-1}\underline{
1^{-1}2^{-1}1^{-1}}1^{-1}2^{-1} \\
&=&1234512341231211^{-1}2^{-1}1^{-1}2^{-1}1^{-1}2^{-1}1^{-1}2^{-1}1^{-1}2^{-1}1^{-1}2^{-1}
\\
&=&123451234123121\left( 1^{-1}2^{-1}\right) ^{6} \\
&=&123451234123121\left( 21\right) ^{-6}.
\end{eqnarray*}

Now using the fact that $t_{c}=\left( 21\right) ^{6}$ we get

\begin{eqnarray*}
t_{b_{0}}t_{b_{1}}t_{b_{2}}t_{c} &=&123451234123121\left(
21\right)
^{-6}\left( 21\right) ^{6}, \\
t_{b_{0}}t_{b_{1}}t_{b_{2}}t_{c} &=&123451234123121,
\end{eqnarray*}

i.e.,
\[
s=t_{b_{0}}t_{b_{1}}t_{b_{2}}t_{c}=t_{c_{1}}t_{c_{2}}t_{c_{3}}t_{c_{4}}t_{c_{5}}t_{c_{1}}t_{c_{2}}t_{c_{3}}t_{c_{4}}t_{c_{1}}t_{c_{2}}t_{c_{3}}t_{c_{1}}t_{c_{2}}t_{c_{1}}.
\]

To see that this is the same expression that we obtained earlier:
\begin{eqnarray*}
s &=&t_{c_{1}}t_{c_{2}}t_{c_{3}}t_{c_{4}}t_{c_{5}}\underline{t_{c_{1}}}%
t_{c_{2}}t_{c_{3}}t_{c_{4}}t_{c_{1}}t_{c_{2}}t_{c_{3}}t_{c_{1}}t_{c_{2}}t_{c_{1}}
\\
 &=&t_{c_{1}}t_{c_{2}}t_{c_{1}}t_{c_{3}}t_{c_{4}}t_{c_{5}}\underline{%
t_{c_{2}}}t_{c_{3}}t_{c_{4}}\underline{t_{c_{1}}}%
t_{c_{2}}t_{c_{3}}t_{c_{1}}t_{c_{2}}t_{c_{1}} \\
 &=&t_{c_{1}}t_{c_{2}}t_{c_{1}}t_{c_{3}}t_{c_{2}}t_{c_{1}}t_{c_{4}}t_{c_{5}}%
\underline{t_{c_{3}}}t_{c_{4}}\underline{t_{c_{2}}}t_{c_{3}}\underline{%
t_{c_{1}}}t_{c_{2}}t_{c_{1}} \\
&=&t_{c_{1}}t_{c_{2}}t_{c_{1}}t_{c_{3}}t_{c_{2}}t_{c_{1}}t_{c_{4}}t_{c_{3}}t_{c_{2}}t_{c_{1}}t_{c_{5}}t_{c_{4}}t_{c_{3}}t_{c_{2}}t_{c_{1}}.
\end{eqnarray*}

For an excellent treatment of this work for general case see
\cite{Ko}.

\section{Main Theorem} \label{mainthm.sec}

\begin{figure}[htbp]
     \centering  \leavevmode
     \psfig{file=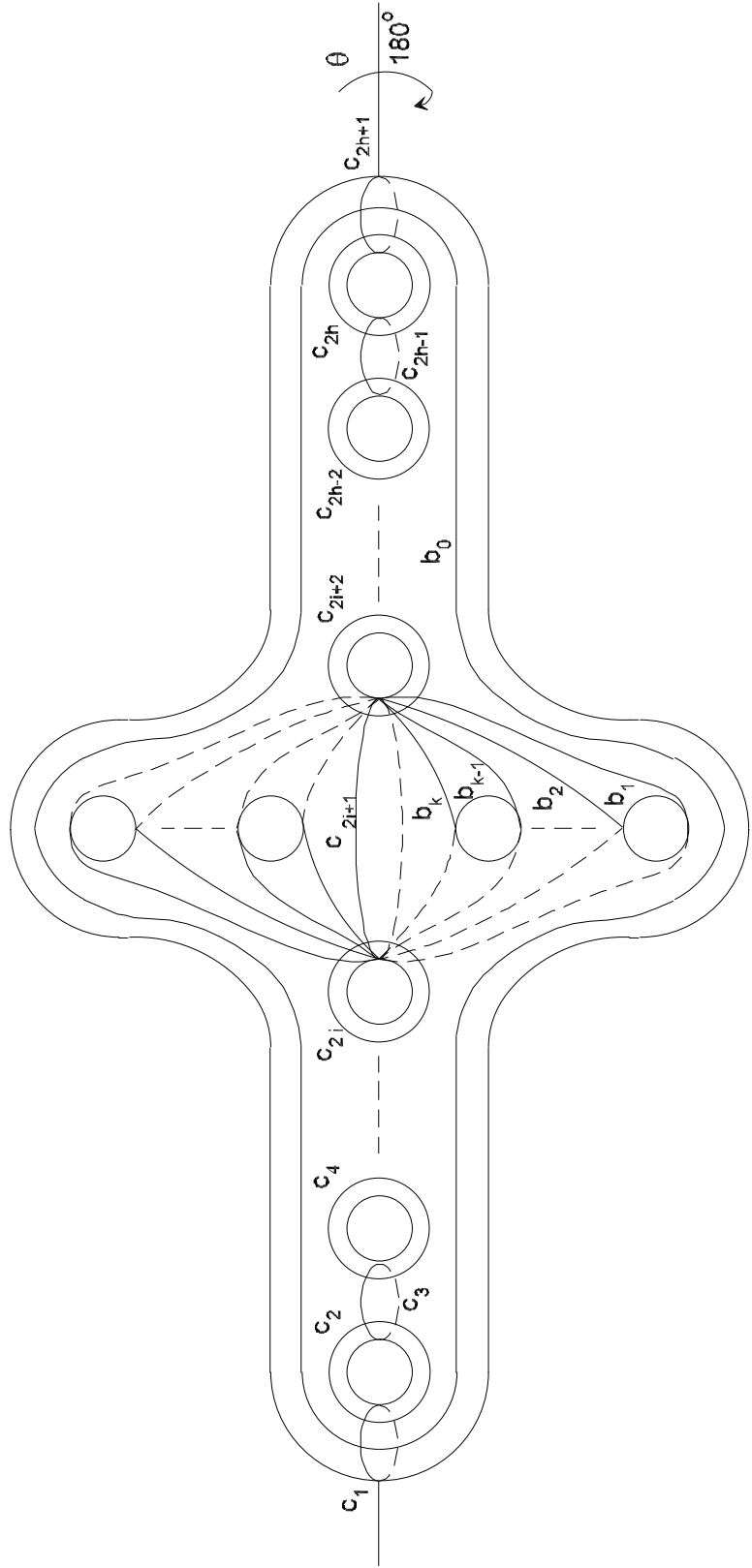,width=14.50in,clip=}
     \caption{The Involution $\theta$}
     \label{theta.fig}
 \end{figure}

The hyperelliptic involution $i:f:\Sigma_{h}\rightarrow
\Sigma_{h}$ and the involution $s:\Sigma_{k}\rightarrow
\Sigma_{k}$ that were described in section \ref{inv.sec}
 can be combined into an
involution $\theta :$ $\Sigma_{h+k}\rightarrow \Sigma_{h+k}$
that supports the action of both $i$ and$\ s$ on bounded subsurfaces of $%
\Sigma_{h+k}$ as shown in Figure \ref{theta.fig}.

\begin{figure}[htbp]
     \centering  \leavevmode
     \psfig{file=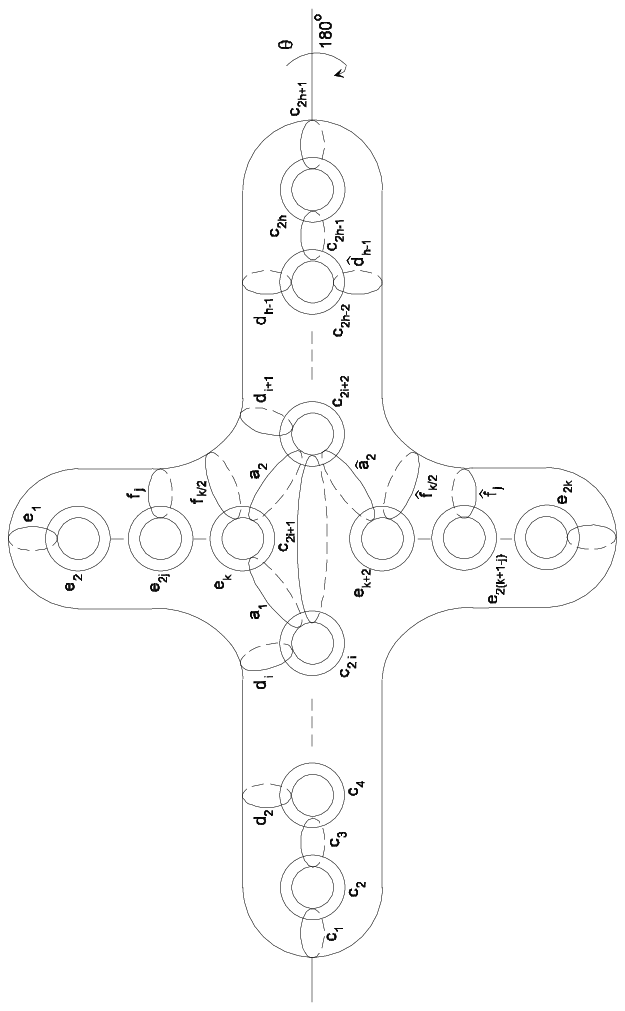,width=14.0in,clip=}
     \caption{Base Cycles}
     \label{basecycles.fig}
 \end{figure}

The main theorem gives a positive Dehn twist expression for
$\theta .$ We will use the same notation for both the cycles and
the Dehn twists about them throughout the section. The order of
the product is from right to left.

\begin{thm} \label{mainthm}
 The positive Dehn twist expression for the
involution that is shown in Figure \ref{theta.fig} is given by
\[
\theta =c_{2i+2}\cdots c_{2h}c_{2h+1}c_{2i}\cdots
c_{2}c_{1}b_{0}c_{2h+1}c_{2h}\cdots c_{2i+2}c_{1}c_{2}\cdots
c_{2i}b_{1}b_{2}\cdots b_{k-1}b_{k}c_{2i+1}.
\]
\end{thm}

\begin{figure}[htbp]
     \centering  \leavevmode
     \psfig{file=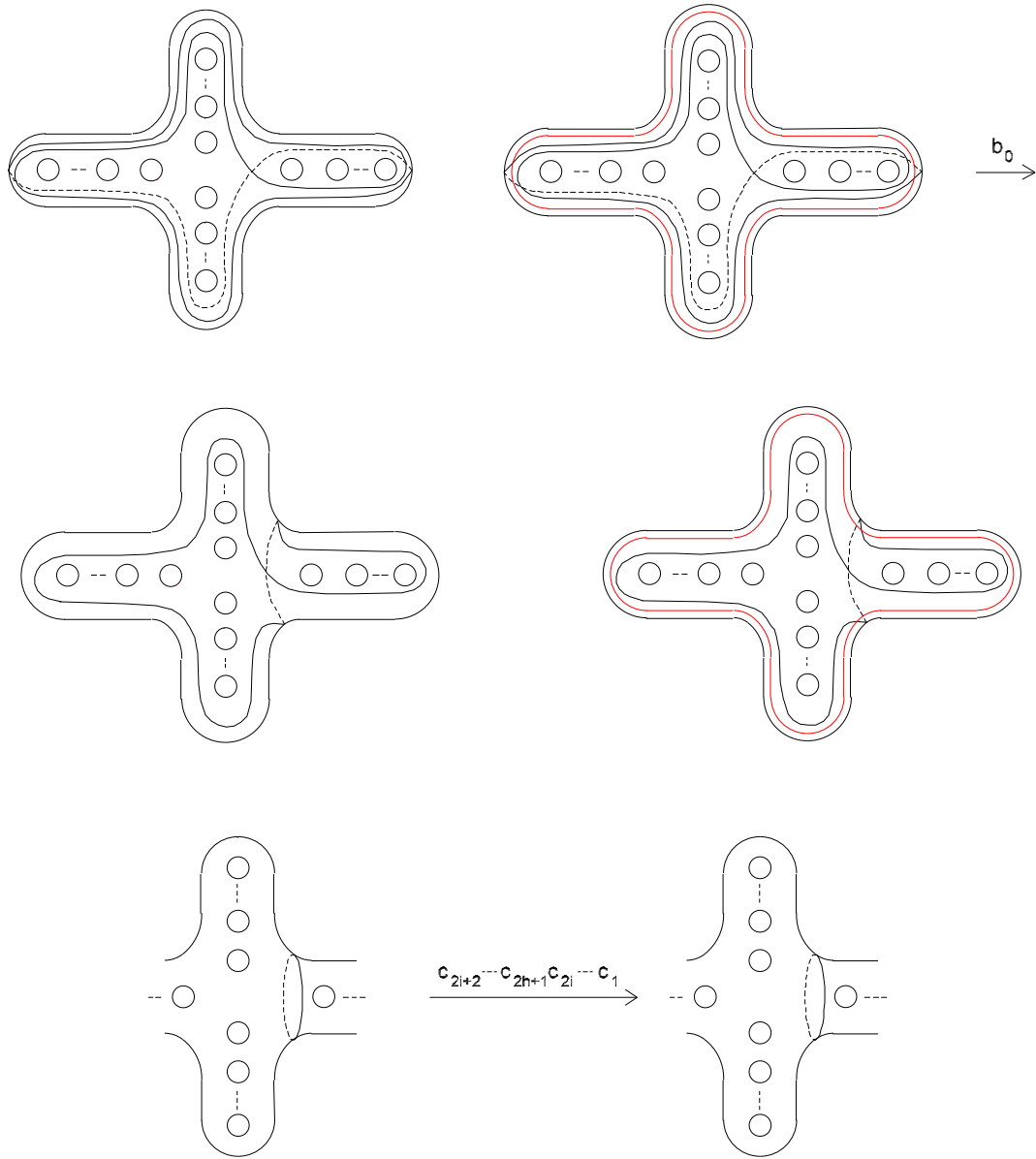,width=12in,clip=}
    \caption{Mapping of a Separating Cycle}
     \label{boundarypage3.fig}
 \end{figure}

\noindent {\bf Proof:} Figure \ref{basecycles.fig} shows a set of
cycles that will constitute a \emph{base  } in the sense that,
mapping of those cycles will ensure the mapping of the subsurfaces
that they bound accordingly. Therefore, it suffices to check the
images of those cycles under $\theta$ and see that they are mapped
to the correct places. Furthermore, for cycle pairs that are
symmetrical with respect to the vertical axis through the center
of the figure, we will look at the mapping of those that are on
the right hand side only. The mapping of their counterpart goes
similarly due to symmetry. Finally, Figure \ref{boundarypage3.fig}
shows an example for the mapping of separating cycles.

\begin{figure}[htbp]
     \centering  \leavevmode
     \psfig{file=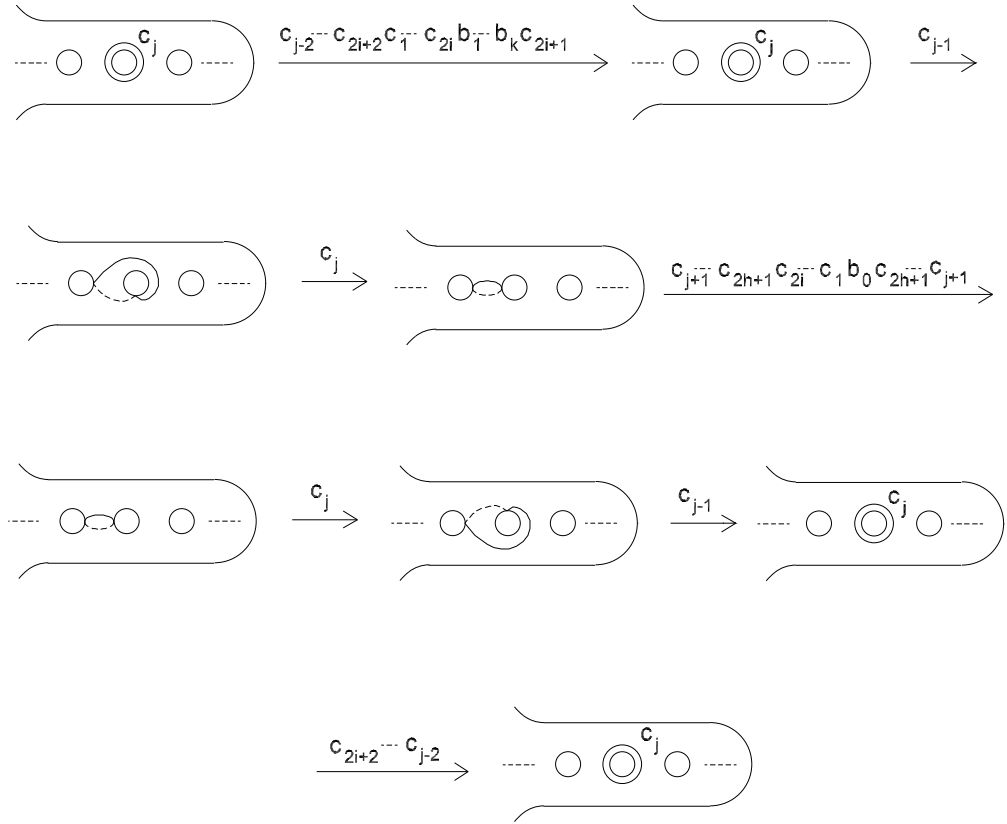,width=8.0in,clip=}
     \caption{Mapping of $c_j$}
     \label{mappingofcj.fig}
 \end{figure}

First we will see the mapping of $c_{j}$ for $j\neq 2i,2i+1$ or
$2i+2.$ Figure \ref{mappingofcj.fig} shows the mapping of such a
cycle for $2i+3\leq j\leq 2h+1$. Whether $c_j$ goes around a hole
or goes through two holes, only two twists are effective in its
mapping: $c_j$ and $c_{j-1}$.\\

\vspace{.3in}

\clearpage

\begin{figure}[htbp]
     \centering  \leavevmode
     \psfig{file=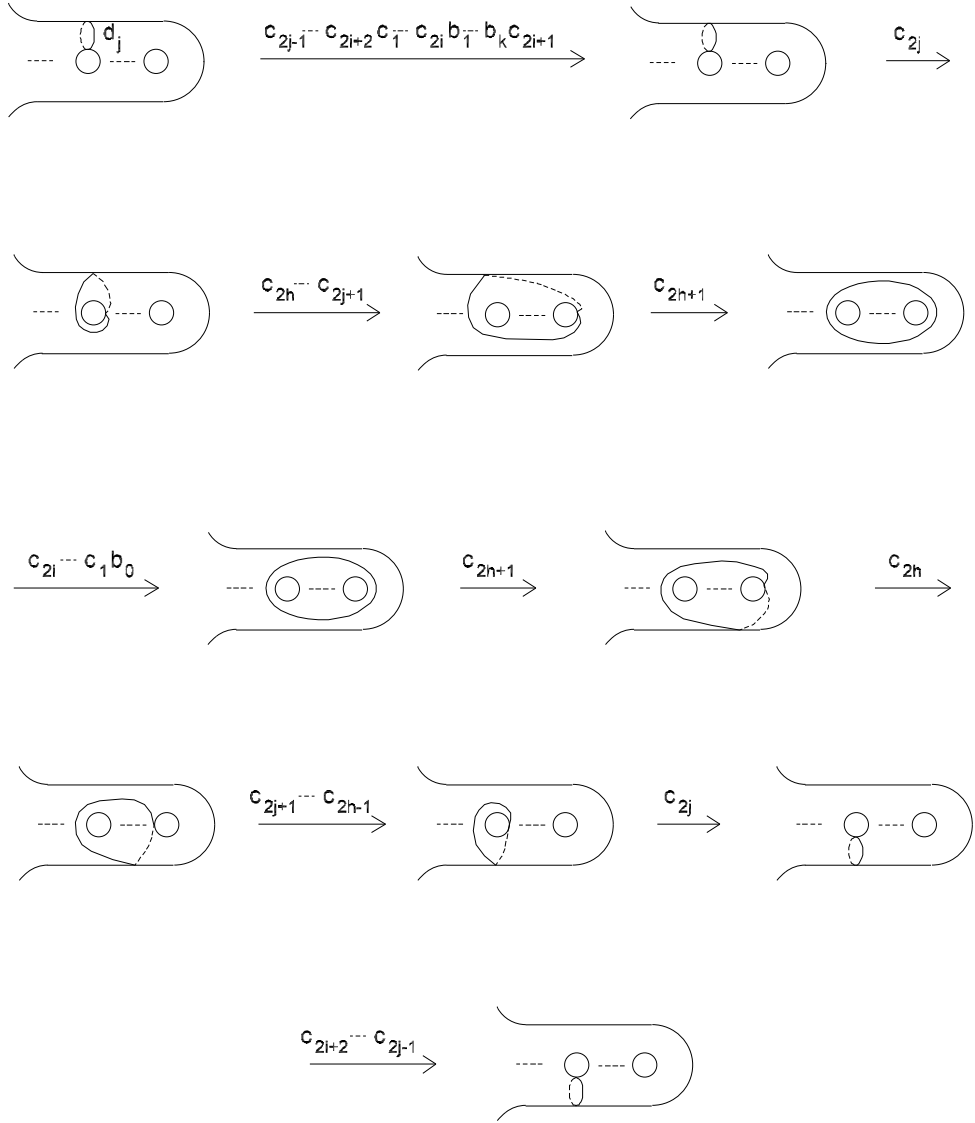,width=11.in,clip=}
     \caption{Mapping of $d_j$}
     \label{mappingofdj.fig}
 \end{figure}

 Next, we will find the image of $d_{j}$ for $2\leq
j\leq h-1.$ Figure \ref{mappingofdj.fig} shows the mapping of
$d_{j}$ for $i+1\leq j\leq h-1.$ As the figure indicates the
twists along $b_j, j=0,1,\ldots,k$ do not take part in mapping of
$d_j$. This is true for $d_i$ and $d_{i+1}$ as well.

\vspace{.75in}

\clearpage

\begin{figure}[htbp]
     \centering  \leavevmode
     \psfig{file=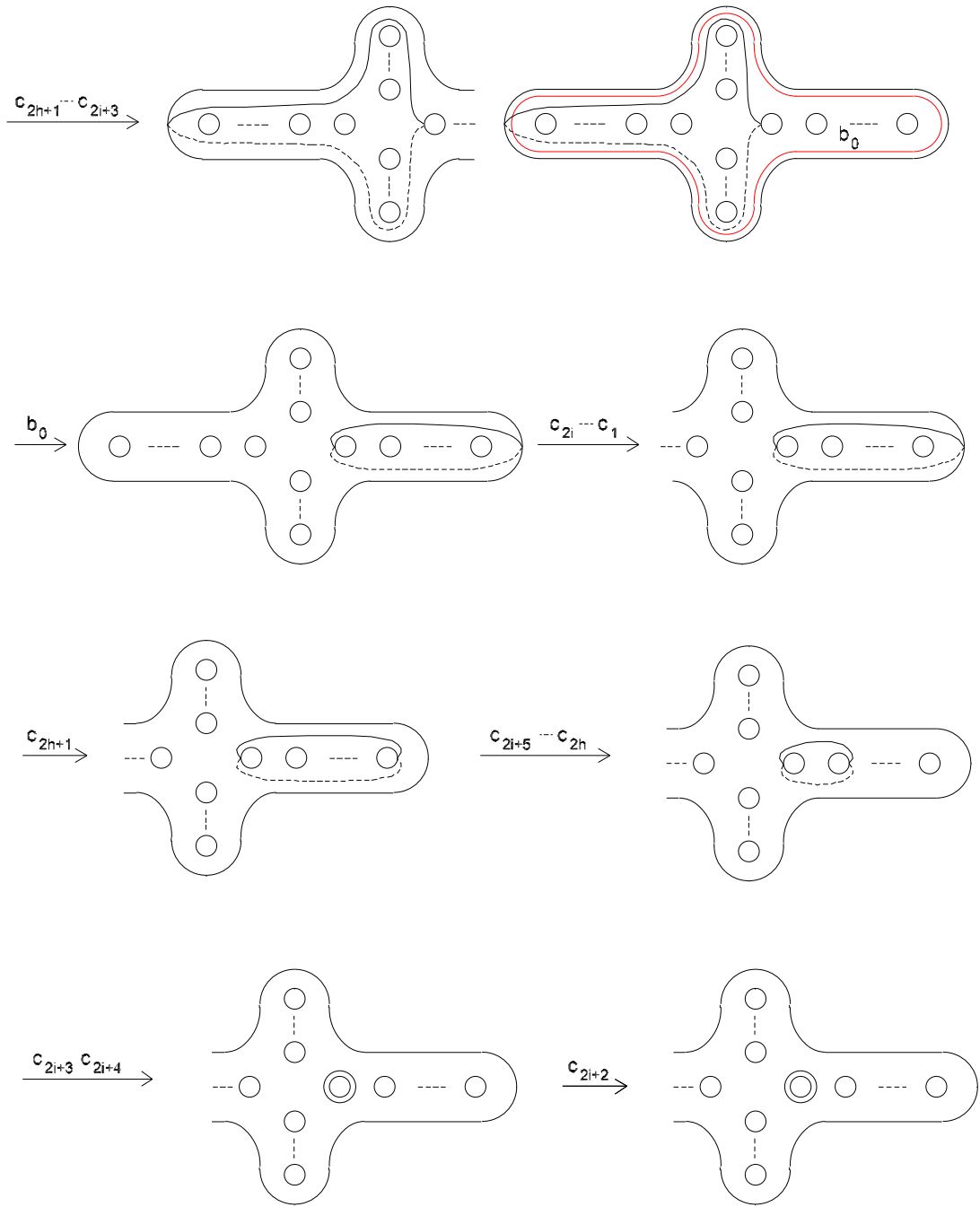,width=13.0in,clip=}
     \caption{Mapping of  $c_{2i+2}$}
     \label{imofc2iplus2page2.fig}
 \end{figure}

Figure \ref{imofc2iplus2page2.fig} shows how $c_{2i+2}$ is mapped.
The mapping of $c_{2i}$ is similar, therefore we omit the proof.
\begin{figure}[htbp]
     \centering  \leavevmode
     \psfig{file=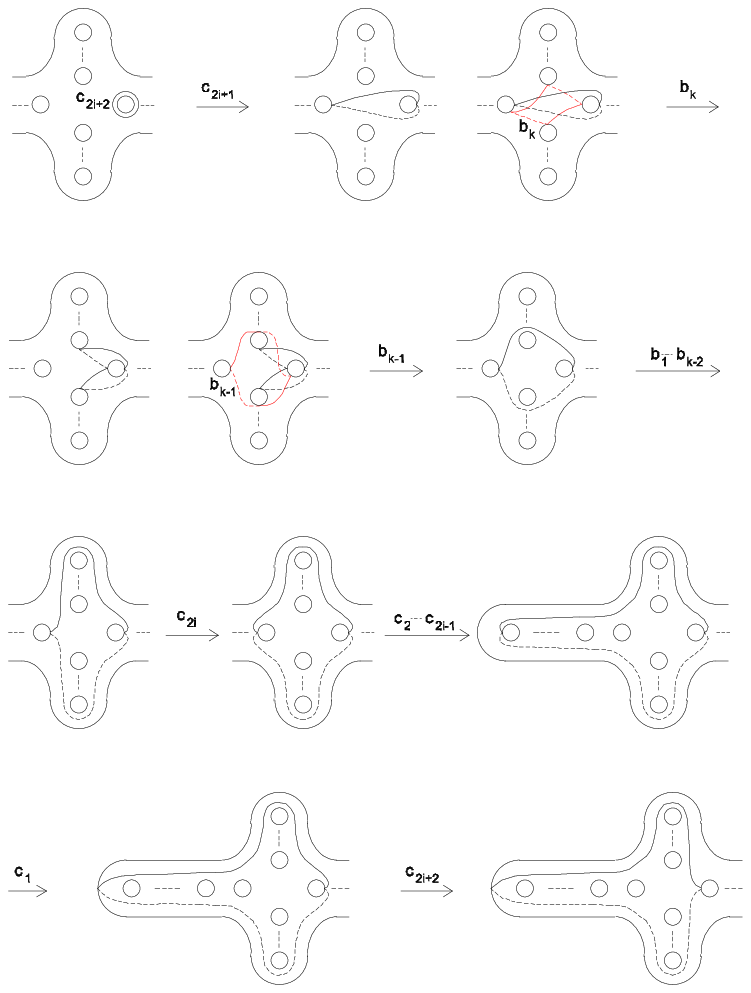,width=13.0in,clip=}
     \label{imofc2iplus2page1.fig}
 \end{figure}

\clearpage

\begin{figure}[htbp]
     \centering  \leavevmode
     \psfig{file=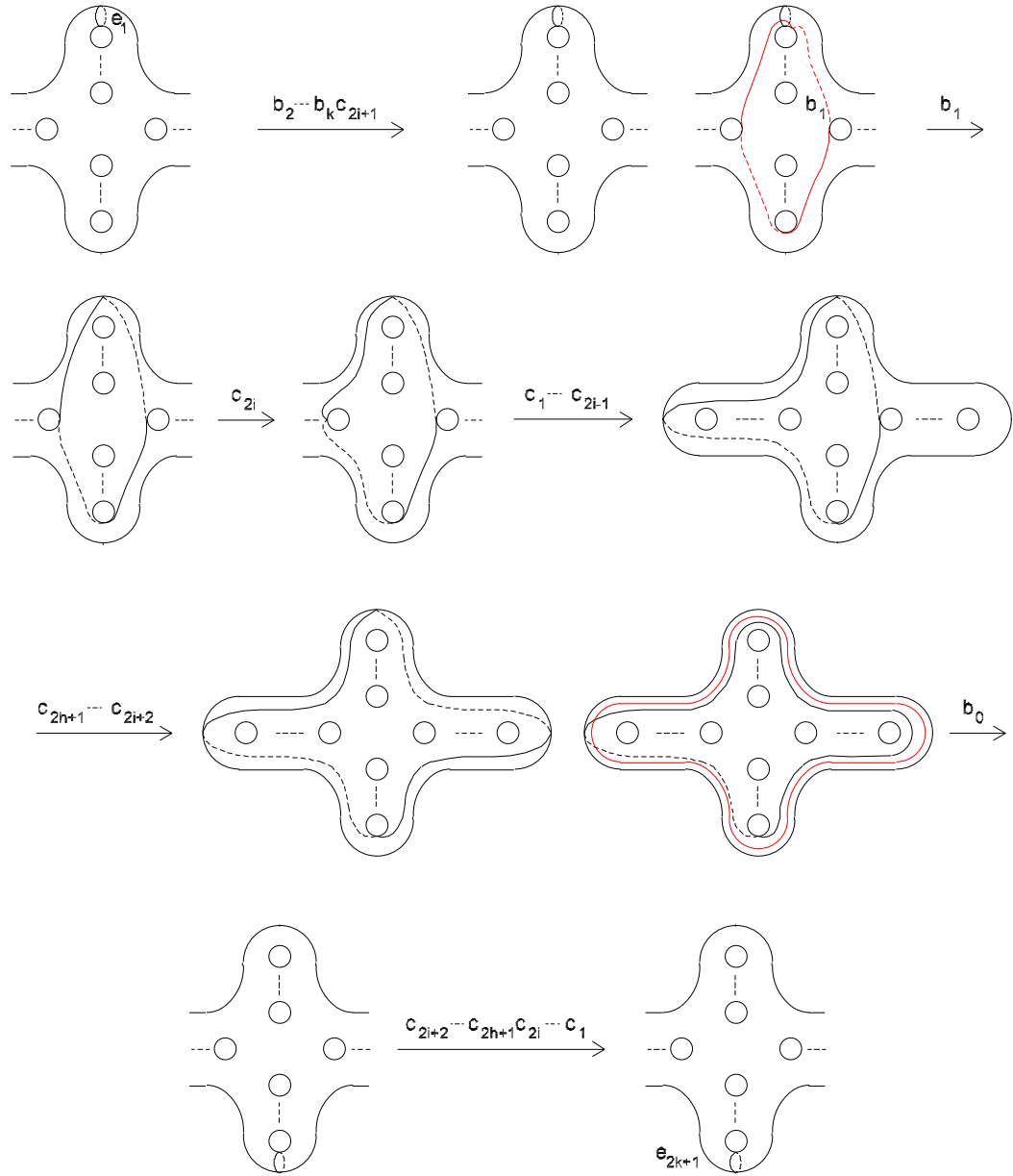,width=12.50in,clip=}
     \caption{Mapping of  $e_{1}$}
     \label{mappingofe1.fig}
 \end{figure}

Figure \ref{mappingofe1.fig} shows the mapping of $e_1.\, b_0$ is
effective in the mapping of $e_1$ as the figure indicates and it
does not take part in the mapping of $e_i, i=2,\ldots,2k. $ The
mapping of $e_{2k+1}$ is the same as that of $e_1$, therefore its
proof is omitted. The cycle to which $b_0$ is applied looks
different in the third line but they are isotopic.

\clearpage

\begin{figure}[htbp]
     \centering  \leavevmode
     \psfig{file=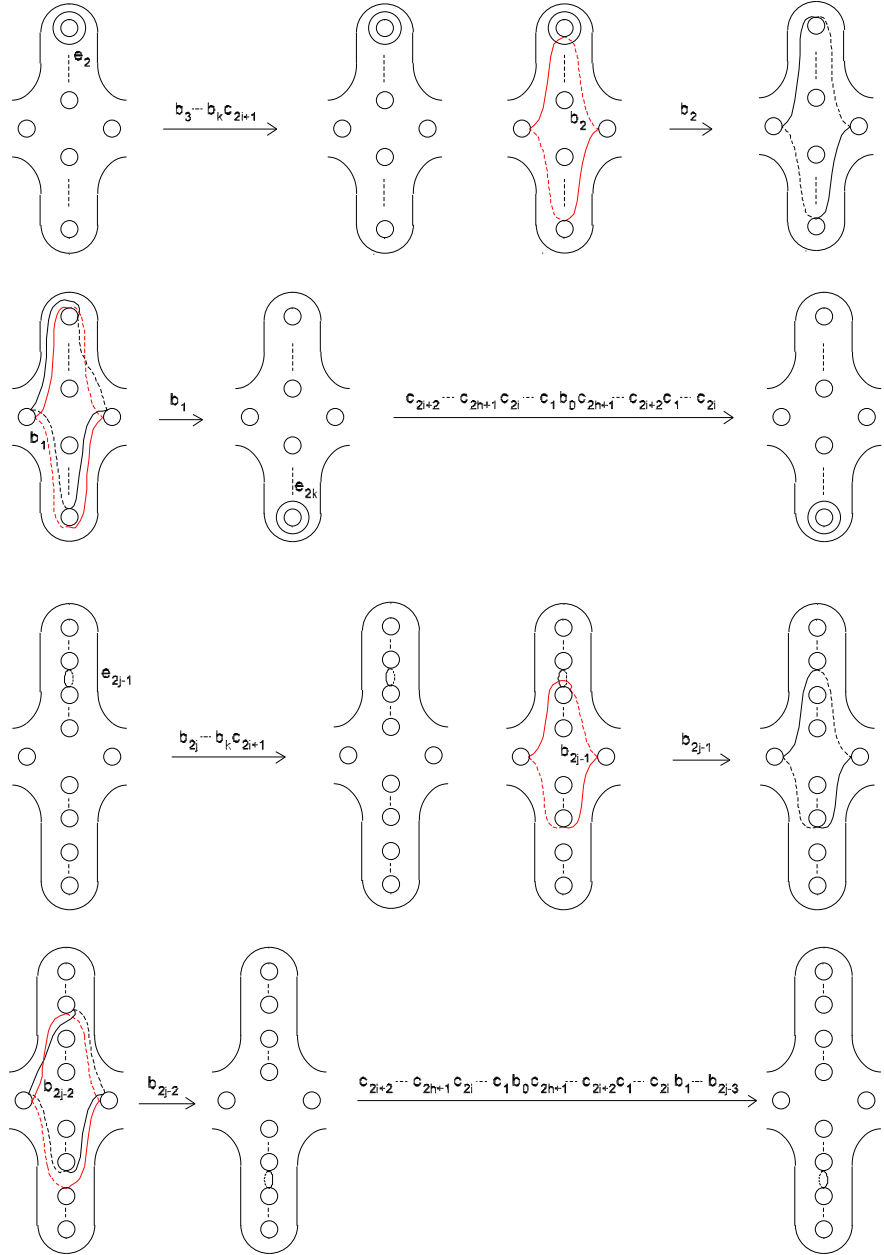,width=12.50in,clip=}
     \caption{Mappings of  $e_2$ and $e_{2j-1},j=2,\ldots,k/2$}
     \label{mappingofe2andej.fig}
 \end{figure}

The mapping of $e_2$ is shown in the first two lines of Figure
\ref{mappingofe2andej.fig}. The mapping of $e_j$ for $j-$ even is
similar to that of $e_2$, including $e_k$. The mapping of
$e_{2j-1},j=2,\ldots,k/2$ is shown in the last two lines of the
same figure. The proof for the mapping of $e_j, j=k+2,\ldots,2k$
is similar, therefore omitted.

\clearpage

\begin{figure}[htbp]
     \centering  \leavevmode
     \psfig{file=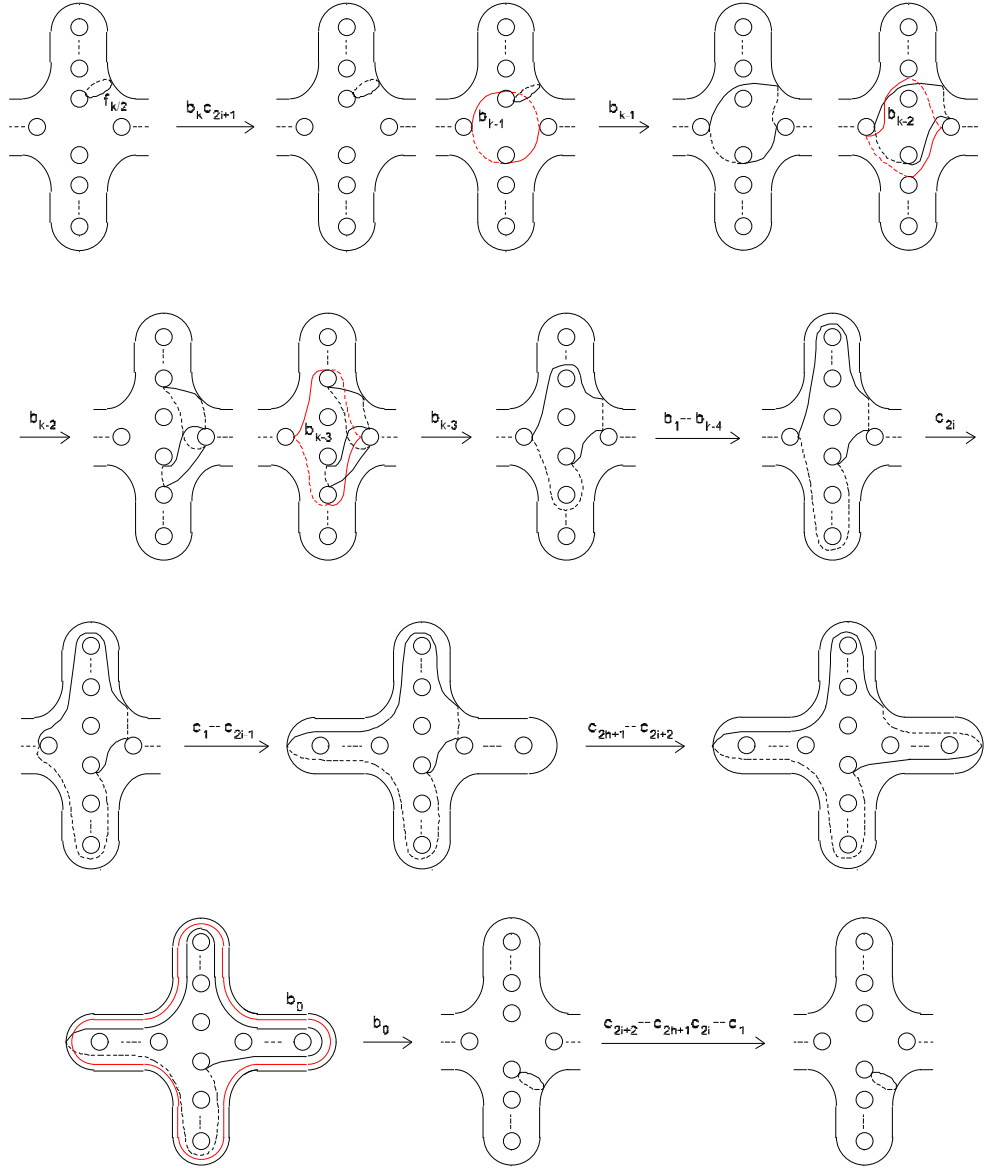,width=15.50in,clip=}
     \caption{Mapping of  $f_{k/2}$}
     \label{mappingoffj.fig}
 \end{figure}

Figure \ref{mappingoffj.fig} shows the mapping of $f_{k/2}$ and
the mapping of $f_j,j=2,\ldots,k/2-1$ are the same. We omit the
proof for $f_j,j=k/2+1,\ldots,k$ due to the same reason. The cycle
to which the twist about $b_0$ is applied in the last line looks
different in the previous figure but they are isotopic.

\clearpage

\begin{figure}[htbp]
     \centering  \leavevmode
     \psfig{file=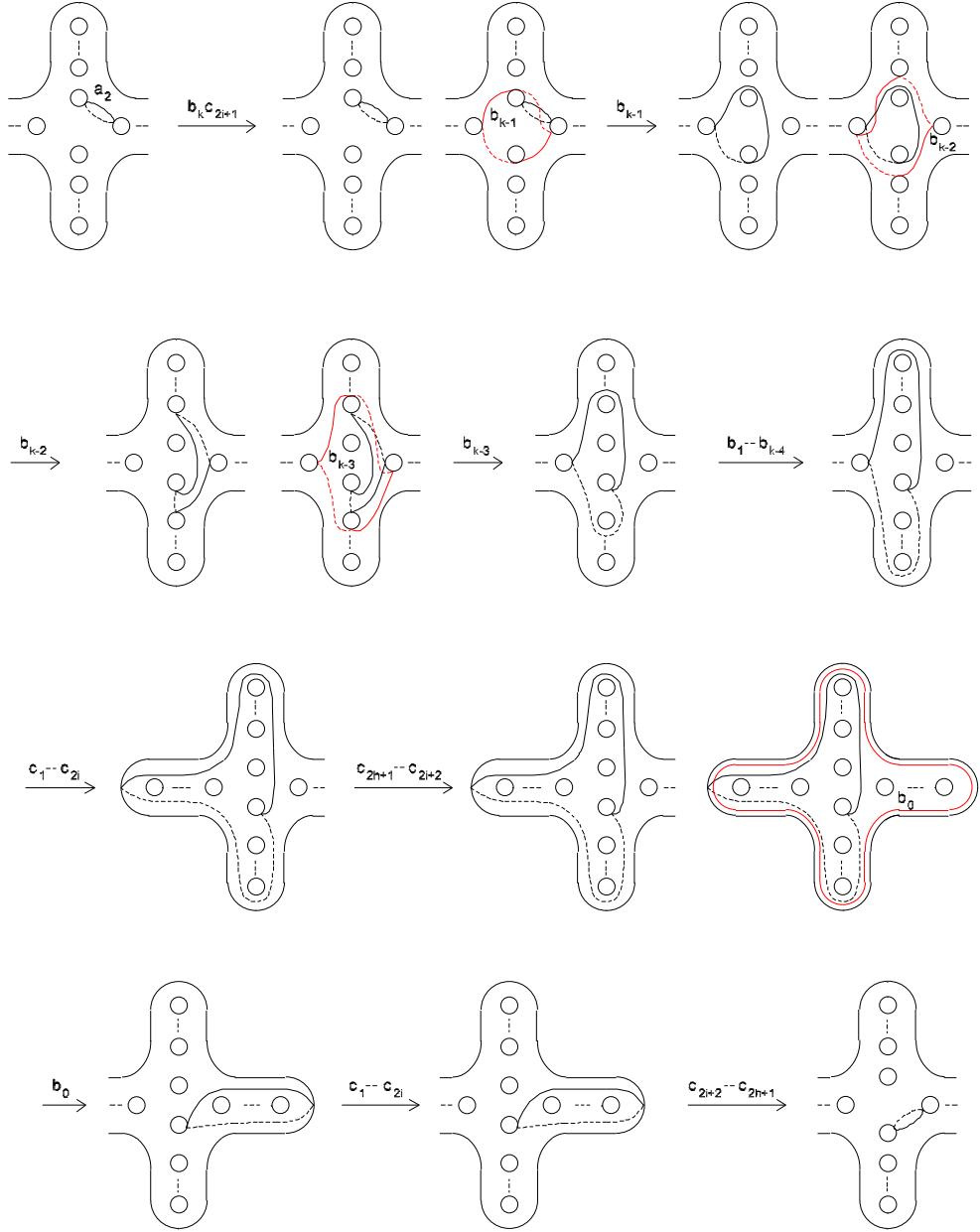,width=13.50in,clip=}
     \caption{Mapping of  $a_{2}$}
     \label{mappingofa2.fig}
 \end{figure}

The mapping of $a_2$ is shown in Figure \ref{mappingofa2.fig}. The
proof for the mapping of its mirror image is omitted because of
the symmetry.

\clearpage

\begin{figure}[htbp]
     \centering  \leavevmode
     \psfig{file=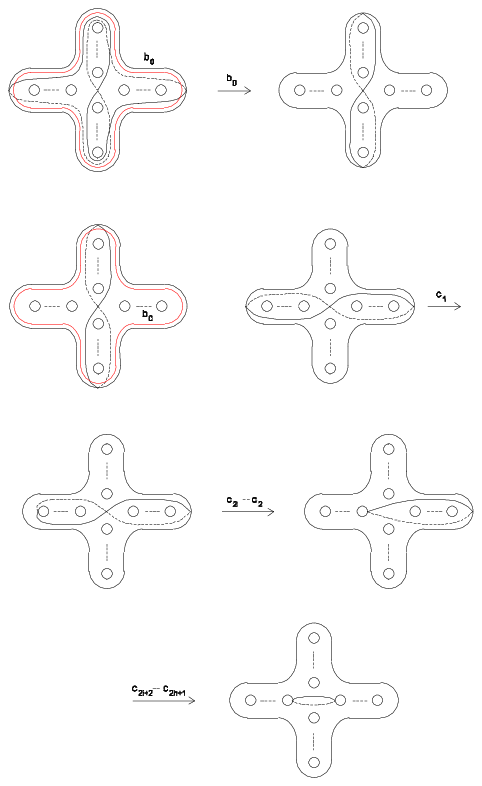,width=17in,clip=}
     \caption{Mapping of  $c_{2i+1}$}
     \label{c2iplus1page2.fig}
 \end{figure}

Figure \ref{c2iplus1page2.fig} shows the mapping of $c_{2i+1}$.
Since  $c_{2i+1} \cap b_k=2$, we have $t_{b_k}(c_{2i+1})=b_k^2\
c_{2i+1}$ by Lemma \ref{intnum.lem}. Also since
$t_{b_k}(c_{2i+1})\ \cap \ b_{k-1}=2$ we have
$t_{b_{k-1}}(t_{b_k}(c_{2i+1}))=b_{k-1}^2\ (b_k^2 \ c_{2i+1}).$
Lemma \ref{intnum.lem} applies to all $b_j,j=0,\ldots , k$ because
each of them intersects the curve they are applied to twice.

\vspace{.2in}
\begin{figure}[htbp]
     \centering  \leavevmode
     \psfig{file=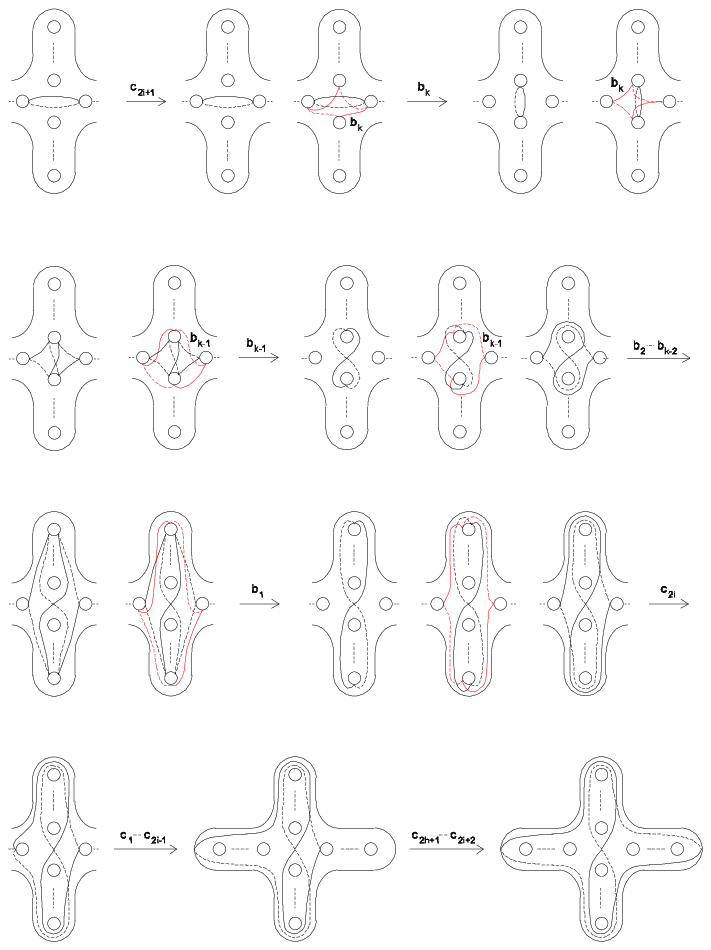,width=13.5in,clip=}
     \label{c2iplus1page1.fig}
 \end{figure}

\clearpage

Figure \ref{boundarypage3.fig} shows the mapping of a separating
cycle.

\vspace{.3in}

\begin{figure}[htbp]
     \centering  \leavevmode
     \psfig{file=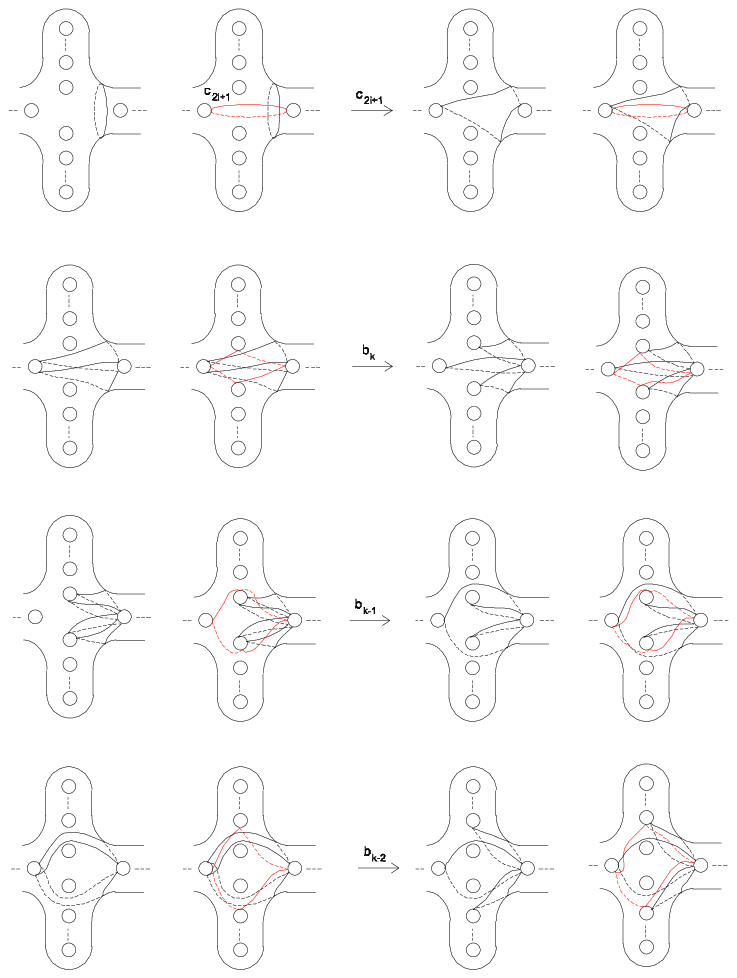,width=13in,clip=}
     \label{boundarypage1.fig}
 \end{figure}

\clearpage

\vspace{1in}

\begin{figure}[htbp]
     \centering  \leavevmode
     \psfig{file=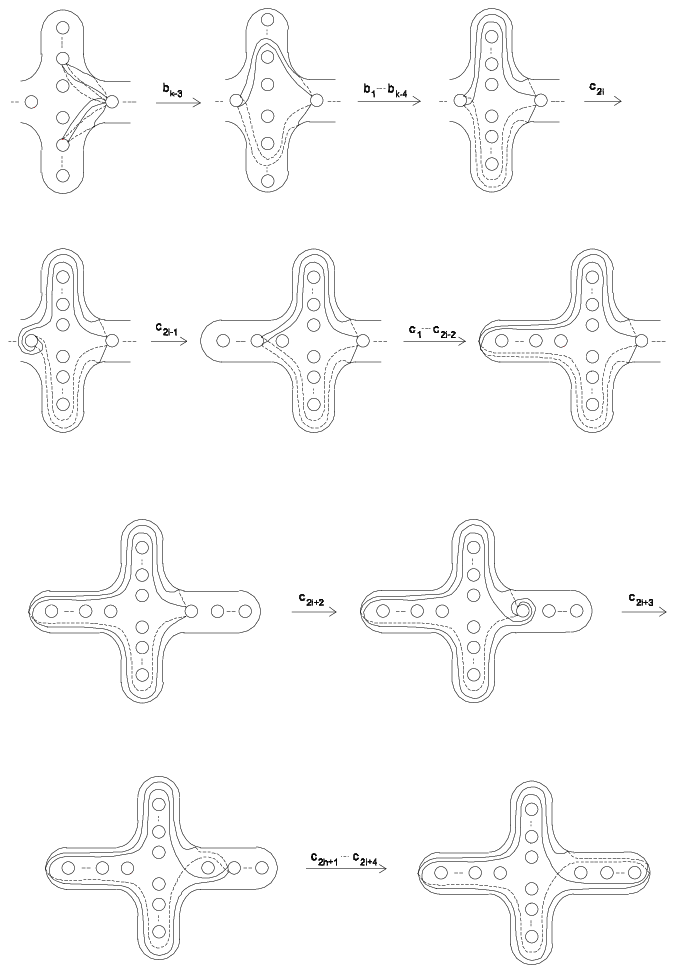,width=15in,clip=}
     \label{boundarypage2.fig}
 \end{figure}

\clearpage

\clearpage

\section{Applications} Next, we will compute the homeomorphism
invariants of the genus $g$ Lefschetz fibrations
\[
 X\longrightarrow S^{2}
\]
 described by the involution $\theta $ that was defined in
Theorem \ref{mainthm}, namely by the word $\theta^2 =1$ in $M_g$,
for some small values of $g$. Consider the surface in Figure
\ref{theta.fig}. Let $k$ denote the genus of the central part of
the figure, which we will call the \emph{vertical genus}. Let
\emph{l} be the total genus on the left and \emph{r} be the total
genus on the right of the vertical component. Let \emph{h} be the
\emph{horizontal genus}, namely the sum of the \emph{left genus}
and the \emph{right genus}; so, $h=l+r.$ If we denote the total
genus by $g$ then $g=h+k=l+r+k$. A quick check reveals that the
total number of twists in the word $\theta^2 =1$ is
$8h+2k+4=2g+4+6h$.

Using the algorithm described in \cite{Oz} we wrote a Matlab
program that computes the signature of the manifold described by
the word $\theta^2 =1$ for given $l,r,$ and $k$. The following
table lists the output of the program for a few values of $l,r,$
and $k$. The first column is for the input of the program in the
format $(l,k,r)$, the second column is for the \emph{horizontal
genus} $h=l+r$, the third column is for the \emph{vertical genus}
$k$, the fourth column is for the total genus $g=h+k$, the fifth
column is for the word length $w=8h+2k+4,$ and the last column is
for the output of the program, which is the signature of the
manifold for the given triple $(l,k,r)$.

\vspace{.3in}

\begin{center}
\begin{tabular}{c|c|c|c|c|c}
(l,k,r) & h=l+r & k & g=h+k & w=8h+2k+4 & signature \\ \hline
 (1,2,1) & 2 & 2 & 4 & 24 &-12 \\
 (1,4,1) & 2 & 4 & 6 & 28 & -12 \\
 (1,6,1) & 2 & 6 & 8 & 32 & -12 \\
 (1,8,1) & 2 & 8 & 10 & 36 & -12 \\
 (1,10,1) & 2 & 10 & 12 & 40 & -12 \\\hline

 (2,2,1) & 3 & 2 & 5 & 32 & -16 \\
 (1,2,2) & 3 & 2 & 5 & 32 & -16 \\
 (2,4,1) & 3 & 4 & 7 & 36 & -16 \\
 (1,4,2) & 3 & 4 & 7 & 36 & -16 \\
 (2,6,1) & 3 & 6 & 9 & 40 & -16 \\
 (1,6,2) & 3 & 6 & 9 & 40 & -16 \\
 (2,8,1) & 3 & 8 & 11 & 44 & -16 \\
 (1,8,2) & 3 & 8 & 11 & 44 & -16 \\ \hline

 (3,2,1) & 4 & 2 & 6 & 40 & -20 \\
 (2,2,2) & 4 & 2 & 6 & 40 & -20 \\
 (1,2,3) & 4 & 2 & 6 & 40 & -20 \\

 (3,4,1) & 4 & 4 & 8 & 44 & -20 \\
 (2,4,2) & 4 & 4 & 8 & 44 & -20 \\
 (1,4,3) & 4 & 4 & 8 & 44 & -20 \\

 (3,6,1) & 4 & 6 & 10 & 48 & -20 \\
 (2,6,2) & 4 & 6 & 10 & 48 & -20 \\
 (1,6,3) & 4 & 6 & 10 & 48 & -20 \\

 (3,8,1) & 4 & 8 & 12 & 52 & -20 \\
 (2,8,2) & 4 & 8 & 12 & 52 & -20 \\
 (1,8,3) & 4 & 8 & 12 & 52 & -20 \\ \hline

\end{tabular}
\end{center}

\clearpage

 The outputs consist of $0$'s and
$-1$'s only. The order in which they appear is the same for a
fixed value of $h$ and $k$ in all the examples above. In other
words the order in which the signature contributions of the $2-$
handles appear does not depend on how the horizontal genus is
distributed to left and right, once we fix $h$ and $k$. Therefore
we will not record $l$ and $r$ for the rest of the examples.

\vspace{.3in}

\begin{center}
\begin{tabular}{c|c|c|c|c}
 h & k & g=h+k & w=8h+2k+4 & signature \\ \hline

 5 & 2 & 7 & 48 &  -24 \\
 5 & 4 & 9 & 52 &  -24 \\
 5 & 6 & 11 & 56 &  -24 \\
 5 & 8 & 13 & 60 &  -24 \\ \hline

 6 & 2 & 8 & 56 &  -28 \\
 6 & 4 & 10 & 60 &  -28 \\
 6 & 6 & 12 & 64 &  -28 \\
 6 & 8 & 14 & 68 &  -28 \\ \hline

 7 & 2 & 9 & 64 &  -32 \\
 7 & 4 & 11 & 68 &  -32 \\
 7 & 6 & 13 & 72 &  -32 \\ \hline

 8 & 2 & 10 & 72 &  -36 \\
 8 & 4 & 12 & 76 &  -36 \\
 8 & 6 & 14 & 80 &  -36 \\ \hline

 \end{tabular}
\end{center}

\vspace{.3in}

 The signature values in the tables above are statistically
convincing that the signature of the Lefschetz fibration given by
the word $\theta^2 =1$  must be
\[
\sigma(X)=-4(h+1).
\]
We do not know a topological proof for that; however, for a proof
of this claim for $h=0,1$ the reader is referred to \cite{Ko}. \\

All of the Dehn twists appearing in the expression for $\theta$
are about nonseparating cycles. Therefore the Lefschetz fibration
that is given by the word $\theta^2=1$ has
\[
8h+2k+4
\]
irreducible singular fibers and its Euler characteristic $\chi(X)$
is

\begin{eqnarray*}
\chi \left( X\right) &=&2\left( 2-2g\right) +8h+2k+4 \\
&=&2\left( 2-2(h+k)\right) +8h+2k+4 \\
&=&8+4h-2k.
\end{eqnarray*}

For the values of the Euler characteristic and the signature
above, $c_1^2(X)$ and $\chi_h(X)$ are

\begin{eqnarray*}
c_1^2(X)&=&3\sigma(X)+2\chi(X) \\
&=&3(-4h-4)+2(8+4h-2k) \\
&=& -4h-4k+4 \\
&=& -4(g-1)
\end{eqnarray*}

and

\begin{eqnarray*}
\chi_h(X)&=&\frac{1}{4}(\sigma(X)+\chi(X)) \\
&=&\frac{1}{4}(8+4h-2k-4h-4) \\
&=& \frac{1}{4}(4-2k)=1-k/2
\end{eqnarray*}

 $\chi_h(X)$ in the above computation makes sense because $X$ has
 almost complex structure. It is an integer since $k$ is even.

 The following are the actual computer outputs for the signature
 computations for the indicated values of $(l,k,r)$.

\vspace{.3in}
\begin{tabular}{c|l}
 $(l,k,r)$ & output \\ \hline
(1,2,1) & 0     +0     +0     +0     +0     +0    -1    -1    -1
-1
-1    -1\\
 & -1    -1     +0    -1    -1    -1    -1     +0     +0     +0     +0
 +0 = -12 \\ \hline
 (1,4,1) & 0     +0     +0     +0     +0     +0     +0     +0    -1    -1    -1    -1
    -1    -1 \\
&  -1    -1     +0     +0     +0    -1    -1    -1    -1     +0 +0
+0     +0     +0  = -12 \\ \hline
 (1,6,1) & 0     +0     +0     +0     +0     +0     +0     +0   +0     +0 -1    -1    -1    -1
    -1    -1 \\
&  -1    -1     +0     +0     +0    +0     +0 -1    -1    -1    -1
+0 +0 +0     +0     +0  = -12 \\ \hline
 (2,2,1) & 0     +0     +0
+0     +0     +0     +0     +0   -1   -1   -1    -1    -1    -1
    -1    -1 \\
&  -1    -1     +0    -1    -1    -1    -1  -1    -1  +0 +0 +0 +0
+0  +0 +0= -16 \\ \hline
 (2,2,2) & 0     +0     +0 +0     +0     +0     +0     +0   +0     +0 -1   -1   -1    -1    -1    -1
   -1    -1  -1    -1 \\
&  -1    -1     +0    -1    -1    -1    -1  -1    -1  -1    -1 +0
+0 +0 +0 +0  +0 +0  +0 +0 \\
& = -20 \\ \hline
 (2,4,2) & 0     +0     +0 +0     +0     +0     +0     +0   +0  +0   +0     +0
 -1   -1   -1    -1    -1    -1    -1    -1  \\
& -1    -1 -1    -1     +0    +0     +0 -1    -1    -1    -1  -1
-1  -1 -1 +0 +0 +0 +0 +0 \\& +0 +0  +0 +0 = -20 \\
\end{tabular}

\begin{tabular}{c|l}
 $(l,k,r)$ & output \\ \hline
 (3,2,2) & 0     +0     +0 +0     +0     +0     +0     +0   +0  +0   +0  +0 -1   -1   -1    -1 \\
 &    -1    -1    -1    -1  -1    -1    -1    -1 -1    -1     +0 -1    -1    -1    -1  -1\\
&-1  -1 -1 -1 -1 +0 +0 +0 +0 +0  +0 +0  +0 +0 +0 +0= -24\\ \hline
 (4,2,4) & 0   +0   +0 +0     +0     +0     +0     +0   +0  +0   +0  +0 +0     +0   +0  +0   +0 +0\\
& -1   -1   -1    -1  -1    -1    -1    -1  -1    -1 -1    -1 -1    -1 -1 -1    -1 -1  \\
 &    -1   -1    +0 -1    -1    -1    -1  -1 -1  -1 -1 -1 -1   -1 -1 -1 -1  -1\\
&-1  +0 +0 +0 +0 +0  +0 +0  +0 +0 +0 +0+0 +0  +0 +0 +0 +0\\
&= -36\\
\hline

\end{tabular}

\vspace{.5in}

{\bf Acknowledgment} I would like to express my deepest gratitude
to my thesis advisor Ronald J. Stern for his continuous support
and guidance.

\clearpage

\clearpage
\bibliographystyle{amsplain}
\pagebreak

\end{document}